\documentclass[a4paper,14pt]{article}
\usepackage[utf8]{inputenc}
\usepackage{amssymb,amsmath,amscd}

\newcommand{\K}{\mathbb{K}}
\newcommand{\Sym}{\mathrm{Sym}}

\newcommand{\Tor}{\mathrm{Tor}}
\newcommand{\odd}{\mathrm{odd}}
\newcommand{\even}{\mathrm{even}}

\newcommand{\Q}{\mathbb{Q}}
\newcommand{\R}{\mathbb{R}}
\newcommand{\Z}{\mathbb{Z}}
\newcommand{\N}{\mathbb{N}}
\newcommand{\Cn}{\mathbb{C}}

\author{DMITRY V. GUGNIN\thanks{This work is supported by the Russian Science Foundation under grant 14-11-00414.}}

\title{\bf \large ON INTEGRAL COHOMOLOGY RING OF SYMMETRIC PRODUCTS}

\date{}

\begin{document}

\pagestyle{plain}
\pagenumbering{arabic}

\maketitle

\begin{abstract}

We prove that the integral cohomology ring modulo torsion $H^*(\Sym^n X;\Z)/\Tor$ for the symmetric product of a connected CW-complex $X$ of finite homology type is a functor of $H^*(X;\Z)/\Tor$ (see Theorem 1). Moreover, we give an explicit description of this functor.

We also consider the important particular case when $X$ is a compact Riemann surface $M^2_g$ of genus $g$. There is a famous theorem of Macdonald of 1962, which gives an explicit description of the integral cohomology ring $H^*(\Sym^n M^2_g;\Z)$. The analysis of the original proof by Macdonald shows that it contains three gaps. All these gaps were filled in by Seroul in 1972, and, therefore, he obtained a complete proof of Macdonald's theorem.  

Nevertheless, in the unstable case $2\le n\le 2g-2$ Macdonald's theorem has a subsection, that needs a slight correction even over $\Q$ (see Theorem 2). 

\end{abstract}

\begin{flushleft}



\end{flushleft}

\medskip

\begin{flushright}
{\it \large To my wife Tanya}
\end{flushright}

\section{Introduction}

The investigation of (co)homology of symmetric products $\Sym^n X:=X^n/S_n$ of CW-complexes $X$ has long history. The first and very deep result is the following theorem of Dold \cite{Dold}.

{\bf Theorem\,(Dold,\,1958).} {\it 
Let $X$ and $Y$ be connected Hausdorff spaces homotopy equivalent to CW-complexes. If $H_*(X;\Z)\cong H_*(Y;\Z)$ (as abelian groups), then $H_*(\Sym^n X;\Z)\cong H_*(\Sym^n Y;\Z)$ for all $n\ge 2.$ Moreover, $H_q(\Sym^n X;\Z), 1\le q < \infty,$ depends only on groups $H_i(X;\Z), 1\le i\le q,$ and the number $n$.}

Evidently, this Dold's theorem implies the following

\textbf{Corollary.} {\it Let $X$ be a connected Hausdorff space homotopy equivalent to CW-complex. If abelian groups $H_q(X;\Z)$ are finitely generated for all $q\ge 1$, then $H_q(\Sym^n X;\Z), 1\le q < \infty$, are also finitely generated for any $n\ge 2$.}

We will say that a connected CW-complex $X$ has finite homology type, if $H_q(X;\Z)$ are finitely generated for all $q\ge 1$.

\textbf{Remark.} The above result of Dold is just a theorem of existence. The rather complicated algorithm for computing $H_q(\Sym^n X;\Z), 1\le q< \infty,$ from $H_i(X;\Z), 1\le i\le q,$ was constructed by Milgram \cite{Mil} in 1969 (for $X$ of finite homology type).

The rational cohomology ring $H^*(\Sym^n X;\Q)$ of connected CW-complexes $X$ of finite homology type can be easily derived from the following classical 

\textbf{Transfer Theorem}. {\it Let $Z$ be an arbitrary simplicial comlex (finite or infinite). Suppose a finite group $G$ acts on $Z$ simplicially, and $\K$ is a field, $\mathrm{char}\K=0$ or $p, (|G|,p)=1$. Then the projection map $\pi:Z\to Z/G$ induces the isomorphism $\pi^*:H^*(Z/G;\K)\cong H^*(Z;\K)^G$.}

The K\"unneth formula for cohomology and simplicial approximation of CW-complexes taken with Transfer Theorem implies the following 

\textbf{Proposition.}  {\it Let $X$ be a connected Hausdorff space homotopy equivalent to a CW-complex of finite homology type. Then the projection map $\pi_n:X^n\to\Sym^n X$ induces the isomorphism 
$$
\pi_n^*:H^*(\Sym^n X;\Q)\cong H^*(X^n;\Q)^{S_n} = (H^*(X;\Q)^{\otimes n})^{S_n}.
$$
}

Let us denote the ring $(H^*(X;\Q)^{\otimes n})^{S_n}$ by $S^n H^*(X;\Q)$. If we know some $\Q$-basis $\gamma_{i}^{q} \in H^q(X;\Q), \forall q\ge 1,$ and rational multiplication table $\gamma_{i}^{q}\gamma_{j}^{r} = c_{q,i;r,j}^{k}\gamma_{k}^{q+r}, c^*_{*}\in \Q$, it is not hard to construct some $\Q$-basis of the $\Q$-algebra $S^n H^*(X;\Q)$ and compute the rational multiplication table for this basis. 

Therefore, the rational cohomology ring $H^*(\Sym^n X;\Q) = S^n H^*(X;\Q)$ of the symmetric product $\Sym^n X$ is a functor of $H^*(X;\Q)$. Moreover, this functor has an explicit form.

{\bf What about the integral cohomology ring $H^*(\Sym^n X;\Z)$? Is it a functor of $H^*(X;\Z)$? And, if it is, how can it be computed?}

It is an open question. Actually, I think the answer is ``NO''. But, even if the ring $H^*(\Sym^n X;\Z)$ is indeed a functor of $H^*(X;\Z)$, it will be a transcendental problem to describe this functor explicitly. The main reason is that even for the very simple case $X=S^4$ the ring $H^*(\Sym^n X;\Z)$ has very complicated torsion when $n$ grows up. Rings with complicated torsion have no good description.

The K\"unneth formula for cohomology gives that for any connected finite CW-complex $X$ the ring $H^*(X^2;\Z)$ has a natural subring $H^*(X;\Z)\otimes H^*(X;\Z)$ and an unnatural direct summand $H_{T}^*(X^2;\Z)$ such that 
$$
H_T^q(X^2;\Z) := \bigoplus_{i+j=q+1} H_{\Tor}^i(X;\Z) \otimes H_{\Tor}^j(X;\Z).
$$
But it is not known what will be the product $\alpha\smile \beta\in H^*(X^2;\Z)$ when $\alpha, \beta \in  H^*_{T}(X^2;\Z)$, or $\alpha \in  H_{T}^*(X^2;\Z)$ and $\beta\in H^*(X;\Z)\otimes H^*(X;\Z)$. Now, we want to pose the following

{\bf Conjecture\,1.} {\it The integral cohomology ring $H^*(X^2;\Z)$ is not a functor of the ring $H^*(X;\Z)$ for connected (even 1-connected) finite CW-complexes $X$. This means that there exist two connected finite CW-complexes $X$ and $Y$ such that $H^*(X;\Z)\cong H^*(Y;\Z)$, but the rings $H^*(X^2;\Z)$ and $H^*(Y^2;\Z)$ are not isomorphic.}

{\bf The main result of this paper is that the ring $H^*(\Sym^n X;\Z)/\Tor$ is a functor of $H^*(X;\Z)/\Tor$. Moreover, we give an explicit description of this functor.} 

Why this result is not a consequence of the classical result over $\Q$? The reason is the following phenomenon.  

{\bf Fact\,$\alpha$.} {\it There exist compact polyhedra (closed manifolds) $L$ and $M$ with torsion-free integral (co)homology having equal rational cohomology rings $H^*(L;\Q)\cong H^*(M;\Q)$, but nonisomorphic integral cohomology rings $H^*(L;\Z)\ncong H^*(M;\Z)$.}

For the simplest example of such manifolds $L$ and $M$ take $S^2\times S^2$ and $\Cn P^2 \# \overline{\Cn P}^2$. Here by $ \overline{\Cn P}^2$ we denote the complex projective plain $\Cn P^2$ equipped with the opposite orientation. 

These four-dimensional 1-connected manifolds have torsion-free integral cohomology, with the multiplicative structure described only by the intersection form. For $S^2\times S^2$ this form is $b_1=\bigl(\begin{smallmatrix}
0& 1\\ 1& 0
\end{smallmatrix}\bigr)$ and for $\Cn P^2 \# \overline{\Cn P}^2$ the form is $b_2=\bigl(\begin{smallmatrix}
1& 0\\ 0& -1
\end{smallmatrix}\bigr)$. The first form $b_1$ is {\it even}, and the second form $b_2$ is {\it odd}. So, they are not isomorphic over $\Z$. Thus, $H^*(S^2\times S^2;\Z)\ncong  H^*(\Cn P^2 \# \overline{\Cn P}^2;\Z)$.

\noindent But, over $\Q$ these forms are isomorphic: 
$$
u=(x,y)\in \Q^2 \Rightarrow b_1(u,u)=2xy=2(\frac{x'+y'}{2})(x'-y')= x'^2 - y'^2.
$$ 
Therefore, we get the isomorphism $H^*(S^2\times S^2;\Q)\cong  H^*(\Cn P^2 \# \overline{\Cn P}^2;\Q)$.  

Here is a series of examples of this phenomenon for compact polyhedra. Fix any $m\ge 1$. Take $S^{2m}$ and mappings $f_{2k}\colon \partial D^{4m} = S^{4m-1} \to S^{2m}$ with Hopf invariant $h(f_{2k}) = 2k, k=1,2,3,\ldots$. Let us take the space $X_{2k}^{4m}:= S^{2m}\bigcup_{f_{2k}} D^{4m}$. The integral cohomology $H^*(X_{2k}^{4m};\Z)$ has the form 
$$
\Z\left<1\right> \oplus \Z \left< u \right>  \oplus \Z \left< v \right>, |u|=2m, |v|=4m, \ \ u^2 = 2kv.
$$
It is evident, that for any $k,l\ge 1, k\neq l$, the integral cohomology rings $H^*(X_{2k}^{4m};\Z)$ and $H^*(X_{2l}^{4m};\Z)$ are not isomorphic. But, for all $k\ge 1,$ the rational cohomology ring $H^*(X_{2k}^{4m};\Q)$ is equal to the following ring 
$$
\Q\left<1\right> \oplus \Q \left< u \right>  \oplus \Q \left< v' \right>, |u|=2m, |v'|=4m, \ \ u^2 = v'.
$$
Of course, for $m=1,2,4$ one may take also the mappings with odd Hopf invariant. 

It is clear, that if $L$ and $M$ are 2-dimensional connected closed orientable manifolds, then the isomorphism $H^*(L;\Q)\cong H^*(M;\Q)$ implies $H^*(L;\Z)\cong H^*(M;\Z)$. 

What about 3-dimensional manifolds? Let $M^3$ be an arbitrary closed connected orientable 3-manifold. Then one has a skew $\Z$-valued 3-form $\mu_{M^3}\colon H^1(M^3;\Z)\times H^1(M^3;\Z) \times H^1(M^3;\Z) \to \Z$ given by the cup-product $\mu_{M^3}\left<a,b,c\right> := abc([M^3])$. The Poincar\'e Duality implies that the ring $H^*(M^3;\Z)/\Tor$ is completely determined by the form $\mu_{M^3}$. The famous theorem of Sullivan \cite{Sul} states that there are no restrictions on the form $\mu_{M^3}$.

\textbf{Theorem\,(Sullivan,\,1975).} {\it Let $H$ be an arbitrary free abelian group of finite rank, and $\mu\colon H\times H\times H \to \Z$ be an arbitrary skew 3-form. Then there is a closed connected orientable 3-manifold $M^3$ such that $H\cong H^1(M^3;\Z)$ and $\mu\cong \mu_{M^3}$.}

Let us take $H:= \Z\left<\alpha_1,\alpha_3,\alpha_3\right>$ and $\mu_s \colon H\times H\times H \to \Z, \  \mu_s\left<\alpha_1,\alpha_3,\alpha_3\right> = s, \  s=1,2,3,\ldots$. By the result of Sullivan there exist 3-manifolds $M^3_s, \ s=1,2,3,\ldots $, with $H\cong H^1(M^3_s;\Z)$ and $\mu_s\cong \mu_{M^3_s}$. Evidently, one may take $M^3_1= T^3$. It is easy to check that the rings $H^*(M^3_s;\Z)/\Tor, s\in \N,$ are pairwise nonisomorphic. But it is also clear that the rational cohomology rings of all these manifolds are the same. Therefore, we have just proved the following

{\bf Fact\,$\beta$.} {\it There exists a series of closed connected orientable 3-manifolds $M_1=T^3, M_2, M_3, \ldots$ having pairwise nonisomorphic integral cohomology rings mod torsion and having one and the same rational cohomology ring.} 

For any connected Hausdorff space $W$ of finite homology type we have a canonical inclusion:
$$
H^*(W;\Z)/\Tor \subset H^*(W;\Q) \  \mbox{and} \ H^*(W;\Q)=(H^*(W;\Z)/\Tor)\otimes \Q.
$$

The most crusual step in the proof of Theorem 1 below is the following Integrality Lemma of Nakaoka \cite{Nak} (see Theorem 2.7 and take $m=1, G=\Z$). It was rediscovered by the author in \cite{my5}. 

{\bf Integrality\,Lemma\,(Nakaoka,\,1957).} {\it Let $X$ be a connected Hausdorff space homotopy equivalent to a CW-complex of finite homology type. Then for any $n\ge 2$ and any ``integral rational'' class $\alpha\in H^*(X;\Z)/\Tor$ the following element $\chi(\alpha) \in S^n H^*(X;\Q)$ lies in the ``integral structure'' subring $H^*(\Sym^n X;\Z)/\Tor:$
$$
\chi(\alpha) = \alpha\otimes 1\otimes\ldots\otimes 1 + 1\otimes\alpha\otimes\ldots\otimes 1 +\ldots + 1\otimes\ldots\otimes 1\otimes\alpha \in H^*(\Sym^n X;\Z)/\Tor.
$$}
The statement of this Lemma is not evident. For instance, if you take any $\beta\in H^{\even}(X;\Z)/\Tor$, which is not of the form $m\gamma$ for $m\ge 2$ and $\gamma\in H^{\even}(X;\Z)/\Tor$, then the following rather simple tensor $\beta\otimes\beta\otimes \ldots\otimes \beta$ is not integral (i.e. it does not belong to $H^*(\Sym^n X;\Z)/\Tor$). This tensor becomes integral only after being multiplied by $n!$.

{\bf In section 2} we prove our main result --- Theorem 1, which gives the explicit description of the ring $H^*(\Sym^n X;\Z)/\Tor$ as a functor of $H^*(X;\Z)/\Tor$. Recently and independently, in the preprint \cite{Ray} Boote and Ray computed the whole integral cohomology ring for symmetric squares of $X=\mathbb{C} P^N$ and $X=\mathbb{H}P^N, N\ge 1$. And they obtained the same answer as we for the torsion free quotient.

Let us denote by $M^2_g$ a compact Riemann surface of genus $g$. {\bf In section 3} we consider a famous Macdonald's theorem of 1962 (see \cite{Mac1}), which gives an explicit description of the integral cohomology ring $H^*(\Sym^n M^2_g;\Z)$. The analysis of the original proof by Macdonald gives that it contains three gaps. All these gaps were filled in by Seroul in 1972 (see \cite{Ser}). Therefore, Seroul in 1972 obtained a complete proof of Macdonald's theorem. 

This correction in Seroul's paper takes 31 pages. But, as we will show the Gap 1 can be filled in by Nakaoka's Integrality Lemma, and the Gap 2 --- by our Theorem 1. The filling in of the Gap 3 takes 7 pages in Seroul's paper. In section 3 we give a more direct 3 pages proof of the statement of the Gap 3. 

Here is a more detailed exposition. At first, Macdonald proves that the ring $H^*(\Sym^n M^2_g;\Z)$ has no torsion. Then he computes the ring $H^*(\Sym^n M^2_g;\Q) = H^*(\Sym^n M^2_g;\Z)\otimes \Q$. Namely, he constructs $2g$ elements $\xi_1,\ldots, \xi_g,\xi'_1,\ldots,\xi'_g\in H^1(\Sym^n M^2_g;\Q)$ and an element $\eta\in H^2(\Sym^n M^2_g;\Q)$ (all of the form $\chi(\alpha)$ for some elements $\alpha\in H^*(M^2_g;\Z)$) such that the ring homomorphism 
$$
f_{\Q}: \Lambda_{\Q} \left< x_1,\ldots,x_g,x'_1,\ldots,x'_g \right>\otimes \Q[y] \to H^*(\Sym^n M^2_g;\Q),
$$
$$
x_i \mapsto \xi_i, \ \  x'_i \mapsto \xi'_i, \ \ 1\le i\le g, \ \ y \mapsto \eta,
$$
is an epimorphism (here, $x_1,\ldots,x_g,x'_1,\ldots,x'_g$ are formal variables of degree $1$, and $y$ is a formal variable of degree $2$). Further Macdonald proves that the ideal $I^*_{Mac, \Q}:= \mathrm{Ker} f_{\Q}$ is generated by the following integral polynomials:
$$
x_{i_1}\ldots x_{i_a}x'_{j_1}\ldots x'_{j_b} (y - x_{k_1}x'_{k_1})\ldots  (y - x_{k_c}x'_{k_c}) y^q, \eqno{(1)}
$$
where $a+b+2c+ q=n+1$ and $i_1,\ldots,i_a,j_1,\ldots,j_b,k_1,\ldots,k_c$ are distinct integers from $1$ to $g$ inclusive.

Then Macdonald implicitly claims that a torsion-free ring can be recovered from its tensor product by $\Q$, which is not always the case. Namely, he claims that the ring homomorphism
$$
f_{\Z}: \Lambda_{\Z} \left< x_1,\ldots,x_g,x'_1,\ldots,x'_g \right>\otimes \Z[y] \to H^*(\Sym^n M^2_g;\Z),
$$
$$
x_i \mapsto \xi_i, \ \  x'_i \mapsto \xi'_i, \ \ 1\le i\le g, \ \ y \mapsto \eta,
$$
is an epimorphism, and the ideal $I^*_{Mac, \Z}:= \mathrm{Ker} f_{\Z}$ is generated by the same polynomials:
$$
x_{i_1}\ldots x_{i_a}x'_{j_1}\ldots x'_{j_b} (y - x_{k_1}x'_{k_1})\ldots  (y - x_{k_c}x'_{k_c}) y^q, 
$$
where $a+b+2c+ q=n+1$ and $i_1,\ldots,i_a,j_1,\ldots,j_b,k_1,\ldots,k_c$ are distinct integers from $1$ to $g$ inclusive.

Here we have three consecutive gaps. 

\textbf{Gap 1.} Why the elements $\xi_1,\ldots, \xi_g,\xi'_1,\ldots,\xi'_g, \eta \in H^*(\Sym^n M^2_g;\Q)$ lie in the integral lattice $H^*(\Sym^n M^2_g;\Z)$? This gap can be filled in by Integrality Lemma (in Macdonald's paper \cite{Mac1} there is no reference to Nakaoka's paper). 

\textbf{Gap 2.} Why the ring homomorphism $f_{\Z}$ is an epimorphism? This gap was filled in by Seroul in \cite{Ser}. It can be also filled in by our Theorem 1 below.

\textbf{Gap 3.} Why the integral polynomials of the form $(1)$ generates the whole ideal $\mathrm{Ker} f_{\Z}$, but not its cocompact sublattice? This gap was filled in by Seroul in \cite{Ser}. In this paper we give a short proof of the statement of the gap 3.

\section{Main Theorem}

\textbf{Theorem 1.} { \it{Let $X$ be a connected Hausdorff space homotopy equivalent to a CW-complex and such that $H_q(X;\Z)$ are finitely generated for all $q\ge 1$.  Denote by $B_m$ the $m$-th rational Betti number of $X$ for all $m\ge 0$. Let $B_{\odd}:=\sum_{m=0}^{\infty} B_{2m+1} \in \Z_{+}\cup \{ \infty \}$ and  $\tilde{B}_{\even}:=\sum_{m=1}^{\infty} B_{2m} \in \Z_{+}\cup \{ \infty \}$. \\
\noindent Take arbitrary additive basises of $H^{2m+1}(X;\Z)/\Tor$ for all $m\ge 0$ and denote by $\alpha_1, \alpha_2,\ldots$ the elements of the union of these basises with the condition $|\alpha_i|\le |\alpha_{i+1}| \ \forall 1\le i \le B_{\odd} - 1$. Also take any additive basises of $H^{2m}(X;\Z)/\Tor$ for all $m\ge 1$ and denote by $\beta_1,\beta_2,\ldots$ the elements of the union of these chosen basises with the same condition $|\beta_i|\le |\beta_{i+1}| \ \forall 1\le i \le \tilde{B}_{\even} - 1$.

\noindent Then for all $n\ge 2$ the following statements hold:

(i) The elements 
$$
\chi(\alpha_i)=\alpha_i\otimes 1\otimes\ldots\otimes 1 + 1\otimes\alpha_i\otimes\ldots\otimes 1 +\ldots + 1\otimes\ldots\otimes 1\otimes\alpha_i \in H^*(\Sym^n X;\Z)/\Tor, 
$$
$$
\chi(\beta_j)=\beta_j\otimes 1\otimes\ldots\otimes 1 + 1\otimes\beta_j\otimes\ldots\otimes 1 +\ldots + 1\otimes\ldots\otimes 1\otimes\beta_j \in H^*(\Sym^n X;\Z)/\Tor, 
$$
$$
1\le i\le B_{\odd}, 1\le j\le \tilde{B}_{\even},
$$
are multiplicative generators of the ring $H^*(\Sym^n X;\Z)/\Tor$. It means that the minimal subring, which contains $\chi(\alpha_i)$ and $\chi(\beta_j)$, coincides with the whole ring $H^*(\Sym^n X;\Z)/\Tor$.

(ii) The elements 
$$
m_1!m_2!\ldots m_l!   \left(   \frac{1}{r! m_1! m_2! \ldots m_l!}   \sum_{  \sigma\in S_n  } \sigma^{-1}(\alpha_{i_1}\otimes \alpha_{i_2}\otimes \ldots \otimes  \alpha_{i_k} \otimes \underbrace{\beta_{j_1} \otimes \ldots \otimes \beta_{j_1}}_{ \mbox{$m_1$ \text{times} } } \otimes \right. 
$$
$$
\left. \underbrace{\beta_{j_2} \otimes \ldots \otimes \beta_{j_2}}_{\mbox{$m_2$ \text{times} } } \otimes \ldots  \otimes \underbrace{ \beta_{j_l} \otimes \ldots \otimes \beta_{j_l} }_{ \mbox{$m_l$ \text{times} } }\otimes \underbrace{  1\otimes \ldots \otimes 1 }_{\mbox{$r$ \text{times} }}   ) \right) \  \in \  H^*(\Sym^n X;\Z)/\Tor, \eqno{(*)} 
$$
where $k\ge 0, l\ge 0, k+l\ge 1, 1\le  i_1<i_2<\ldots <i_k \le B_{\odd}, 1\le j_1<j_2<\ldots <j_l \le \tilde{B}_{\even}, m_1\ge 1,m_2\ge 1,\ldots, m_l\ge 1,  r\ge 0, k+m_1+m_2+\ldots + m_l + r = n$,

\noindent form an additive basis of $H^*(\Sym^n X;\Z)/\Tor$. The multiplier $\frac{1}{r! m_1! m_2! \ldots m_l!}$ inside big brackets kills the repetitions of elementary tensors under the summation symbol.

Moreover, if we know the integral multiplication table for the basis $\{ \alpha_i,\beta_j \}_{1 \le i\le B_{\odd}, 1\le j \le \tilde{B}_{\even}}$ of $H^*(X;\Z)/\Tor$, then there is a simple algorithm for computing the integral multiplication table for the considered basis $(*)$ of the ring $H^*(\Sym^n X;\Z)/\Tor$. } }

\textbf{Corollary 1.} {\it Let $X$ and $Y$ be connected Hausdorff spaces homotopy equivalent to CW-complexes and having finitely generated integral homology groups in all dimensions. If the rings $H^*(X;\Z)/\Tor$ and $H^*(Y;\Z)/\Tor$ are isomorphic, then the rings $H^*(\Sym^n X;\Z)/\Tor$ and $H^*(\Sym^n Y;\Z)/\Tor$ are also isomorphic for all $n\ge 2$.}

\textbf{Proof.}

We will use the following notations:
$$
A^*_{\Z}:= H^*(X;\Z)/\Tor,  \  A^*:=A^*_{\Z}\otimes \Q=H^*(X;\Q), 
$$
$$ 
D^*_{n\Z}:=H^*(\Sym^n X;\Z)/\Tor,  \  D^*_n:=D^*_{n\Z}\otimes \Q =H^*(\Sym^n X;\Q),
$$
$$
S^nA^*:= (A^{\otimes n})^{S_n}= (H^*(X^n;\Q))^{S_n}.
$$

One has canonical inclusions $A^*_{\Z}\subset A^*$ and $D^*_{n\Z}\subset D^*_n$. The projection $\pi_n:X^n\rightarrow \Sym^n X$ induces the canonical isomorphism $\pi_n^*:H^*(\Sym^n X;\Q)=D^*_n\rightarrow S^nA^*=(H^*(X^n;\Q))^{S_n}$. We will identify the rings $D^*_n$ and $S^nA^*$ via this isomorphism $\pi^*_n$. 

Let us write an explicit formula for the (right) action of the symmetric group $S_n$ on $A^{\otimes n}$.  For any $k,l\in \mathbb{N}$ we set $\varepsilon_{kl}=1$ if $k>l$, and $\varepsilon_{kl}=0$ if $k\le l$. Then for any homogeneous $a_1,a_2,\ldots, a_n\in A^*$ one has:
$$
(a_1\otimes a_2\otimes \ldots\otimes a_n)\sigma  \  = \    \sigma^{-1} (a_1\otimes a_2\otimes \ldots\otimes a_n) \  =   
$$
$$
= (-1)^{\sum_{1\le p< q\le n} \varepsilon_{\sigma^{-1}(p)\sigma^{-1}(q)} |a_p||a_q| } a_{\sigma(1)}\otimes a_{\sigma(2)}\otimes \ldots\otimes a_{\sigma(n)}   \  \  \forall \sigma\in S_n.
$$

First, let us prove that the elements
$$
\chi(\alpha_{i_1},\alpha_{i_2},\ldots,\alpha_{i_k}  |  \underbrace{\beta_{j_1},\ldots,\beta_{j_1}}_{\mbox{$m_1$ \text{times}}},\ldots, \underbrace{\beta_{j_l},\ldots,\beta_{j_l}}_{\mbox{$m_l$ \text{times}}}  ) :=
$$
$$
= \frac{1}{(n-(k+m_1+\ldots+m_l))!} \sum_{\sigma\in S_n} \sigma^{-1} (\alpha_{i_1} \otimes \ldots \otimes  \alpha_{i_k} \otimes \underbrace{\beta_{j_1} \otimes \ldots \otimes \beta_{j_1}}_{ \mbox{$m_1$ \text{times} } } \otimes \ldots  \otimes \underbrace{ \beta_{j_l} \otimes \ldots \otimes \beta_{j_l} }_{ \mbox{$m_l$ \text{times} } }\otimes  1\otimes \ldots \otimes 1 )
$$
of the type $(*)$ lie in the subring $D^*_{n\Z}\subset D^*_n=S^nA^*$. For this it is sufficient to prove that for any homogeneous $\xi_1,\ldots,\xi_n\in A^{\odd}_{\Z}$ and $\eta_1,\ldots,\eta_n\in A^{\even \ge 2}_{\Z}$ the elements
$$
\chi(\xi_1,\ldots,\xi_k | \eta_1,\ldots,\eta_s) := \frac{1}{(n-k-s)!} \sum_{\sigma \in S_n}  \sigma^{-1}(\xi_1\otimes \ldots\otimes \xi_k \otimes \eta_1\otimes \ldots\otimes \eta_s\otimes 1\otimes \ldots \otimes 1)  \  \in \  D^*_n,
$$
where $0\le k,s\le n, 1\le k+s\le n,$ lie in the subring $D^*_{n\Z}$. For the proof we will use the finite induction on pairs $(k,s)$. 

The base case
$$
\chi(\xi_1| \varnothing) := \frac{1}{(n-1)!} \sum_{\sigma \in S_n}  \sigma^{-1}(\xi_1\otimes 1\otimes \ldots \otimes 1)  =  \xi_1\otimes 1\otimes \ldots \otimes 1 + \ldots +  1\otimes \ldots \otimes 1\otimes  \xi_1 = \chi(\xi_1) \in   D^*_{n\Z},
$$ 
$$
\chi(\varnothing | \eta_1 ) := \frac{1}{(n-1)!} \sum_{\sigma \in S_n}  \sigma^{-1}(\eta_1\otimes 1\otimes \ldots \otimes 1)  =  \eta_1\otimes 1\otimes \ldots \otimes 1 + \ldots +  1\otimes \ldots \otimes 1\otimes  \eta_1 = \chi(\eta_1)  \in   D^*_{n\Z},
$$
follows from Integrality Lemma. 

The inductive step: suppose that for all pairs $(k,s), 1\le k+s\le q< n,$ the statement holds true. Let us show that for all pairs $(k,s), k+s =q+1,$ the statement also holds. Clearly, it is sufficient to derive the statement for $(k+1,s)$ and $(k,s+1)$, where $k+s=q$, from the statements for all pairs $(k,s), k+s=q$.

\textbf{Case} $(k+1,s)$. Suppose $k,s\ge 1$. One has the following computation:
$$
\chi(\xi_1,\ldots,\xi_k | \eta_1,\ldots,\eta_s) \chi(\xi_{k+1} | \varnothing)    = 
$$
$$
=    \frac{1}{(n-k-s)!} \frac{1}{(n-1)!} \sum_{\sigma,\tau \in S_n} \sigma^{-1}(\xi_1\otimes \ldots\otimes \xi_k\otimes \eta_1\otimes\ldots\otimes \eta_s\otimes 1\otimes \ldots\otimes 1) \tau^{-1}(\xi_{k+1}\otimes 1\otimes \ldots\otimes 1)    =
$$
$$
=   \frac{1}{(n-k-s)!} \frac{1}{(n-1)!} \sum_{\sigma,\tau \in S_n} \sigma^{-1}(\xi_1\otimes \ldots\otimes \xi_k\otimes \eta_1\otimes\ldots\otimes \eta_s\otimes 1\otimes \ldots\otimes 1)  (\tau\sigma)^{-1}  (\xi_{k+1}\otimes 1\otimes \ldots\otimes 1)   =
$$
$$
=   \frac{1}{(n-k-s)!} \frac{1}{(n-1)!} \sum_{\sigma,\tau \in S_n} \sigma^{-1}( (\xi_1\otimes \ldots\otimes \xi_k\otimes \eta_1\otimes\ldots\otimes \eta_s\otimes 1\otimes \ldots\otimes 1)  \tau^{-1} ( \xi_{k+1}\otimes 1\otimes \ldots\otimes 1) )  =
$$
$$
=   \frac{1}{(n-k-s)!}  \sum_{\sigma \in S_n} \sigma^{-1}( (\xi_1\otimes \ldots\otimes \xi_k\otimes \eta_1\otimes\ldots\otimes \eta_s\otimes 1\otimes \ldots\otimes 1)  ( \xi_{k+1}\otimes 1\otimes \ldots\otimes 1  + \ldots + 1\otimes \ldots \otimes 1\otimes \xi_{k+1}  ) )   =
$$
$$
=   \frac{1}{(n-k-s)!} \sum_{\sigma \in S_n} \left(  \phantom{1^{1^1}}  (-1)^{k-1}  \sigma^{-1}( \xi_1 \xi_{k+1}\otimes\xi_2\otimes  \ldots\otimes \xi_k\otimes \eta_1\otimes\ldots\otimes \eta_s\otimes 1\otimes \ldots\otimes 1) \ + \right.
$$
$$
+ \   (-1)^{k-2}  \sigma^{-1}( \xi_1 \otimes\xi_2\xi_{k+1}\otimes \xi_{3}\otimes  \ldots\otimes \xi_k\otimes \eta_1\otimes\ldots\otimes \eta_s\otimes 1\otimes \ldots\otimes 1) \ + \ \ldots \ +
$$
$$
+ \  (-1)^0  \sigma^{-1}( \xi_1 \otimes \xi_2\otimes  \ldots\otimes \xi_k\xi_{k+1} \otimes \eta_1\otimes\ldots\otimes \eta_s\otimes 1\otimes \ldots\otimes 1)  \  +
$$
$$
+ \  \sigma^{-1}( \xi_1\otimes  \ldots\otimes \xi_k \otimes \eta_1\xi_{k+1} \otimes \eta_2\otimes  \ldots\otimes \eta_s\otimes 1\otimes \ldots\otimes 1)  \  + \  \ldots \  + 
$$
$$
+ \  \sigma^{-1}( \xi_1\otimes  \ldots\otimes \xi_k \otimes \eta_1 \otimes \eta_2\otimes  \ldots\otimes \eta_s\xi_{k+1} \otimes 1\otimes \ldots\otimes 1)  \   +
$$
$$
\left.   + \ (n-k-s)\sigma^{-1}( \xi_1\otimes  \ldots\otimes \xi_k \otimes \eta_1 \otimes \eta_2\otimes  \ldots\otimes \eta_s\otimes \xi_{k+1} \otimes 1\otimes \ldots\otimes 1) \phantom{1^{1^1}}   \right)  =
$$
$$
=  (-1)^{k-1} \chi(\xi_2,\xi_3,\ldots,\xi_k |  \xi_1\xi_{k+1}, \eta_1, \eta_2,\ldots,\eta_s) + (-1)^{k-2} \chi(\xi_1,\xi_3,\ldots,\xi_k |  \xi_2\xi_{k+1}, \eta_1, \eta_2,\ldots,\eta_s) + 
$$
$$
+\ldots + (-1)^0 \chi(\xi_1,\xi_2,\ldots,\xi_{k-1} |  \xi_k\xi_{k+1}, \eta_1, \eta_2,\ldots,\eta_s) + 
$$
$$
+  \chi(\xi_1,\xi_2,\ldots,\xi_k, \eta_1\xi_{k+1} |  \eta_2, \eta_3,\ldots,\eta_s) + \ldots + \chi(\xi_1,\xi_2,\ldots,\xi_k, \eta_s\xi_{k+1} |  \eta_1, \eta_2,\ldots,\eta_{s-1}) + 
$$
$$
+ \chi(\xi_1,\xi_2,\ldots,\xi_k,\xi_{k+1} |  \eta_1, \eta_2,\ldots,\eta_s).
$$
It is easy to see that the above computation implies the desired integrality property $ \chi(\xi_1,\xi_2,\ldots,\xi_k,\xi_{k+1} |  \eta_1, \eta_2,\ldots,\eta_s)  \in D^*_{n\Z}$. 

Let $k=0, 1\le s\le n-1$. By means of the similar computation one has the following identity, which proves the integrality property:
$$
\chi(\varnothing | \eta_1,\ldots,\eta_s)\chi(\xi_1 | \varnothing) =  \chi(\eta_1\xi_1 |  \eta_2, \eta_3,\ldots,\eta_s) + \ldots + \chi(\eta_s\xi_1 |  \eta_1, \eta_2,\ldots,\eta_{s-1})  +  \chi(\xi_1 | \eta_1,\eta_2,\ldots, \eta_s).
$$

Suppose $1\le k\le n-1, s=0$. Then one has:
$$
\chi(\xi_1,\xi_2,\ldots,\xi_k | \varnothing) \chi(\xi_{k+1} | \varnothing)  =  (-1)^{k-1}\chi(\xi_2,\xi_3,\ldots,\xi_k | \xi_1\xi_{k+1}) +  
$$
$$
+ (-1)^{k-2}\chi(\xi_1,\xi_3,\ldots,\xi_k | \xi_2\xi_{k+1})  + \ldots + (-1)^0 \chi(\xi_1,\xi_2,\ldots,\xi_{k-1} | \xi_k\xi_{k+1}) + \chi(\xi_1,\xi_2,\ldots,\xi_k,\xi_{k+1} | \varnothing ). 
$$

Therefore, in the current case the inductive step is proved.

\textbf{Case} $(k,s+1)$. Using the similar computation one can derive the following formula, which proves the inductive step in this case: 
$$
\chi(\xi_1,\ldots,\xi_k | \eta_1,\ldots,\eta_s) \chi( \varnothing | \eta_{s+1})    = 
$$
$$
=  \chi(\xi_1 \eta_{s+1},\xi_2,\ldots,\xi_k  |  \eta_1,\ldots,\eta_s) + \chi(\xi_1 ,\xi_2\eta_{s+1},\xi_3,\ldots,\xi_k  |  \eta_1,\ldots,\eta_s)  + \ldots +
$$
$$
+  \chi(\xi_1,\xi_2,\ldots,\xi_k\eta_{s+1}  |  \eta_1,\ldots,\eta_s) + \chi(\xi_1,\xi_2,\ldots,\xi_k  |  \eta_1\eta_{s+1},\eta_2,\ldots,\eta_s) + \ldots +
$$
$$
+  \chi(\xi_1,\xi_2,\ldots,\xi_k  |  \eta_1,\eta_2,\ldots,\eta_s\eta_{s+1}) +  \chi(\xi_1,\xi_2,\ldots,\xi_k  |  \eta_1,\eta_2,\ldots,\eta_s,\eta_{s+1}).
$$

Therefore, we have just proved that the elements of the type $(*)$ lie in the subring $D^*_{n\Z}\subset D^*_n$.

It follows from the above induction reasoning, that the minimal subring $\Z[\chi(\alpha_i),\chi(\beta_j)]\subset D^*_{n\Z}$, generated by $\chi(\alpha_i),\chi(\beta_j), 1\le i\le B_{\odd}, 1\le j\le \tilde{B}_{\even}$, contains all elements of the type $(*)$. Thus, the statement $(ii)$ of the theorem implies the statement $(i)$. Let us prove the part $(ii)$. 

It is evident, that the graded $\Q$-vector space $S^nA^* = D^*_n$ could be linearly generated by the elements $\Sym(\omega) = \frac{1}{n!}\sum_{\sigma\in S_n}\sigma^{-1}(\omega)$, where $\omega = a_1\otimes a_2\otimes \ldots\otimes a_m\otimes 1\otimes\ldots \otimes 1$ and $a_1,\ldots,a_m$ are homogeneous elements of $A^{*\ge 1}, 1\le m\le n$. As for odd-dimensional elements $a_i$ we have $a_i=p_i^s \alpha_s$, and for even-dimensional elements $a_i$ we have $a_i=q_i^t \beta_t$, where $p^*_*,q^*_*\in \Q$, therefore $S^nA^*$ is linearly generated by the elements
$$
\sum_{  \sigma\in S_n  } \sigma^{-1}(\alpha_{i_1}\otimes \alpha_{i_2}\otimes \ldots \otimes  \alpha_{i_k} \otimes \underbrace{\beta_{j_1} \otimes \ldots \otimes \beta_{j_1}}_{ \mbox{$m_1$ \text{times} } }  \otimes \ldots  \otimes \underbrace{ \beta_{j_l} \otimes \ldots \otimes \beta_{j_l} }_{ \mbox{$m_l$ \text{times} } }\otimes \underbrace{  1\otimes \ldots \otimes 1 }_{\mbox{$r$ \text{times} }}   ), 
$$
where $k\ge 0, l\ge 0, k+l\ge 1, 1\le  i_1\le i_2\le \ldots \le i_k \le B_{\odd}, 1\le j_1<j_2<\ldots <j_l \le \tilde{B}_{\even}, m_1\ge 1,m_2\ge 1,\ldots, m_l\ge 1,  r\ge 0, k+m_1+m_2+\ldots + m_l + r = n$.\\
But, if $i_s=i_{s+1}$ for some $1\le s\le k-1$, then the corresponding symmetric element is equal to zero. Thus, we have proved that the elements of the type $(*)$ linearly generates $S^nA^*=D^*_n$.

Moreover, it can be derived from Macdonald's formula for Betti numbers of symmetric products (see \cite{Mac2}), that these elements $(*)$ are $\Q$-linear independent. But, we will show it independently, and at the same time we will prove, that the elements $(*)$ are $\Z$-basis of the ``integral structure'' subring $D^*_{n\Z}\subset D^*_n$.

One has $\Z$-basis 
$$
\alpha_1,\alpha_2,\ldots \ \in H^1(X;\Z)/\Tor \oplus H^3(X;\Z)/\Tor \oplus \ldots.
$$
Let 
$$
a^1,a^2,\ldots \ \in H_1(X;\Z)/\Tor \oplus H_3(X;\Z)/\Tor \oplus \ldots
$$
be the dual $\Z$-basis: $\left< \alpha_i, a^j \right> = \alpha_i(a^j) = \delta^j_i$.

Also one has $\Z$-basis 
$$
\beta_1,\beta_2,\ldots \ \in H^2(X;\Z)/\Tor \oplus H^4(X;\Z)/\Tor \oplus \ldots.
$$
Let 
$$
b^1,b^2,\ldots \ \in H_2(X;\Z)/\Tor \oplus H_4(X;\Z)/\Tor \oplus \ldots
$$
be the dual $\Z$-basis: $\left< \beta_i, b^j \right> = \beta_i(b^j) = \delta^j_i$.

Consider the following element:
\begin{multline*}
c = \frac{1}{m'_1! m'_2! \ldots m'_{l'}!}  a^{i'_1}\otimes a^{i'_2} \otimes \ldots \otimes  a^{i'_{k'}} \otimes \\
\otimes\underbrace{ b^{j'_1} \otimes \ldots \otimes b^{j'_1} }_{ \mbox{$m'_1$ \text{times} } } \otimes  \ldots  \otimes \underbrace{ b^{j'_{l'}} \otimes \ldots \otimes b^{j'_{l'}} }_{ \mbox{$m'_{l'}$ \text{times} } }\otimes  \underbrace{  1\otimes \ldots \otimes 1 }_{\mbox{$r'$ \text{times} }}  \in H_*(X^n;\Q),   \ \ \ \ \ \ (2)
\end{multline*} 
where $k'\ge 0, l'\ge 0, k'+l'\ge 1, 1\le  i'_1 < i'_2 < \ldots < i'_{k'} \le B_{\odd}, 1\le j'_1<j'_2<\ldots <j'_{l'} \le \tilde{B}_{\even}, m'_1\ge 1,m'_2\ge 1,\ldots, m'_{l'}\ge 1,  r'\ge 0, k'+m'_1+m'_2+\ldots + m'_{l'} + r' = n$.

Let $\bar{c}=\pi_{n*}(c) \in H_*(\Sym^n X;\Q)$, where $\pi_n:X^n \to \Sym^n X$ is a canonical projection. We will prove now, that the element $\bar{c}$ lies in ``integral structure'' $H_*(\Sym^n X;\Z)/\Tor \subset H_*(\Sym^n X;\Q)$. This fact could be easily derived from the following lemma.

\textbf{Lemma 1.}{\it{ Let $X$ be a connected Hausdorff space homotopy equivalent to a CW-complex and such that $H_q(X;\Z)$ are finitely generated for all $q\ge 1$.  Suppose $b\in H_{2p}(X;\Z)/\Tor, p\ge 1.$ Then for all $m\ge 2$ the element $\frac{1}{m!}b\otimes \ldots \otimes b\in H_*(X^m;\Q)$ is mapped onto ``integral'' element $\pi_{m*}(\frac{1}{m!}b\otimes \ldots \otimes b)\in H_*(\Sym^m X;\Z)/\Tor\subset H_*(\Sym^m X;\Q)$, where $\pi_m:X^m\to \Sym^m X$ is a canonical projection.}}

\textbf{Remark}. This lemma is a weaker version of the so called {\it divided power operations} for integral homology of symmetric products. These operations are mentioned in the  paper \cite{Dold2} by Dold (see pp. 15-16).

\textbf{Proof of lemma 1.} It is well known, that any integral homology class $a\in H_q(X;\Z), q\ge 1,$ can be represented as an image $a = f_*([M^q])$ of a fundamental class $[M^q]$ of some connected oriented topologically normal $q$-dimensional pseudomanifold $M^q$, where $f:M^q\to X$ is an appropriate map. 

A pseudomanifold $M^q$, which is always assumed to be a compact simplicial complex, is called topologically normal, if the subcomplex $Z(M^q)\subset M^q$ of singular points of $M^q$ does not locally disconnect the pseudomanifold $M^q$ (the point $x\in M^q$ is called singular if there is no neibourhood $x\in U_x$ which is homeomorphic to $\R^q$). If $q=1,2$ then any topologically normal pseudomanifold $M^q$ is a topological manifold. For all $q\ge 3$ there exist topologically normal pseudomanifolds $M^q$ with singular points. 

It is easy to check that for all $m\ge 2$ the symmetric product $\Sym^m M^1$ of a 1-dimensional manifold $M^1$ is an $m$-dimensional topological manifold with boundary. The fundamental natural homeomorphism $\Sym ^m \mathbb{C}\cong \mathbb{C}^m, m\ge 2,$ implies the fact that symmetric products $\Sym^m M^2$ of 2-dimensional topological manifolds are indeed $2m$-dimensional topological (even smoothable) manifolds. If $M^2$ is closed and oriented, then $\Sym^m M^2$ is also closed and oriented. 

If $q\ge 3$, then symmetric products $\Sym^m \R^q, m\ge 2,$ all have singular points. Moreover, these spaces are not even homology manifolds (there exist points $x\in \Sym^m \R^q$, in which local homology is not the same as for $\R^{mq}$). But, it is not hard to check, that symmetric products $\Sym^m M^q, m\ge 2, q\ge 3,$ of connected topologically normal $q$-dimensional pseudomanifolds $M^q$ are also connected topologically normal pseudomanifolds of dimension $mq$. It is also easy to prove, that pseudomanifold $\Sym^m M^q, m\ge 2, q\ge 3,$ is orientable {\it iff} $q$ is even and $M^q$ is orientable. 

Suppose $b\in H_{2p}(X;\Z)/\Tor, p\ge 1.$ Then there exists a connected oriented topologically normal $2p$-dimensional pseudomanifold $M^{2p}$ and a map $f:M^{2p}\to X$ such that $b=f_*([M^{2p}])$. We have the following commutative diagram:
$$
\begin{CD}
M^{2p}\times\ldots\times M^{2p}  @>f\times\ldots\times f>>  X \times \ldots \times X \\
@VV\pi_m^{M}V   @VV\pi_m^XV  \\
\Sym^m M^{2p} @>\Sym^m(f)>> \Sym^m X
\end{CD}
$$
It is evident, that $[M^{2p}\times\ldots\times M^{2p}] = [M^{2p}]\otimes\ldots\otimes [M^{2p}]$ and $\pi_{m*}^M([M^{2p}\times\ldots\times M^{2p}]) = m! [\Sym^m M^{2p}]$. 
The above commutative diagram implies the equality
$$
\pi_{m*}^X (\frac{1}{m!} b\otimes \ldots \otimes b)  =  \pi_{m*}^X (\frac{1}{m!} f_*\otimes\ldots\otimes f_* ([M^{2p}\times\ldots\times M^{2p}])) =
$$
$$
= \Sym^m(f)_* ( \pi_{m*}^M  ( \frac{1}{m!} [M^{2p}\times\ldots\times M^{2p}])  ) = \Sym^m(f)_* ( [\Sym^m M^{2p}] ) \in H_{2pm}(\Sym^m X;\Q),
$$
which concludes the proof of the lemma. \  \  \ \ \ \ \ \ \ \ \ \ \ \ \ \ \ \ \ \ \ \ \ \ \ \ \ \ \ \ \ \ \ \ \ \ \ \ \ \ \ \ \ \ \ \ \ \ \ \ \ \ \ \ \ \ \ \ \ \ \ \ \ \ \ \ \ \ \ \ \ \ \ \ \ \ \ \ \ \ \ \ \ \ \ \ \ \ \ \ \ \ \ \ \ \ \ \ \ \ \ \ \ \ \ \ \ $\Box$

Therefore, we have the above considered element $c\in H_*(X^n;\Q)$ with the ``integral'' image $\bar{c}=\pi_{n*}(c)\in H_*(\Sym^n X;\Z)/\Tor$. Let $\bar{\gamma}\in H^*(\Sym^n X;\Z)/\Tor$ be an arbitrary element of the type $(*)$. Denote also $\gamma=\pi_{n}^*(\bar{\gamma})\in H^*(X^n;\Z)/\Tor$. We have $\left< \bar{\gamma}, \bar{c} \right> = \left< \gamma, c \right>\in \Z$.

Now, let $\bar{\gamma}\in H^*(\Sym^n X;\Z)/\Tor$ runs over all elements of type $(*)$, and $c\in H_*(X^n;\Q)$ runs over all elements of the type $(2)$. It is easy to see, that the value $\left< \bar{\gamma}, \bar{c} \right> = \left< \gamma, c \right>$ always belongs to the set $\{0,+1,-1\}$. Moreover, $\left< \bar{\gamma}, \bar{c} \right> \ne 0$ {\it iff} the elements $\gamma$ and $c$ are constructed from the same combinatorial data, i.e. 
$$
k'=k, i'_1=i_1,\ldots,i'_k=i'_{k'}, l'=l, m'_1=m_1,\ldots,m'_{l'}=m_l, j'_1=j_1, \ldots, j'_{l'}=j_l, r'=r.
$$ 

Let us fix any $p\ge 1$. Then, by construction the elements $\bar{\gamma} \in H^*(\Sym^n X;\Z)/\Tor$ of degree $p$ forms a $\Q$-basis of $H^p(\Sym^n X;\Q)$ and the elements $\bar{c}\in H_*(\Sym^n X;\Z)/\Tor$ of degree $p$ forms a $\pm$-dual $\Q$-basis of $H_p(\Sym^n X;\Q)$. But, as we take the elements $\bar{\gamma}$ and $\bar{c}$ from ``integral'' lattices, this fact implies that our elements of the type $(*)$ with degree $p$ forms a $\Z$-basis of the lattice $H^p(\Sym^nX;\Z)/\Tor$. 

The statement $(ii)$ of the theorem is proved. As was mentioned above, $(ii)$ implies $(i)$, so the part $(i)$ is also proved. By product, we have shown that the elements of the type $(2)$ with degree $p$ forms a $\Z$-basis of the homology lattice $H_p(\Sym^nX;\Z)/\Tor$, i.e. we have constructed also an explicit $\Z$-basis in homology modulo torsion.

Suppose we know the integral multiplication table for the considered $\Z$-basis $\alpha_i,\beta_j,1\le i\le B_{\mathrm{odd}}, 1\le j \le \tilde{B}_{\mathrm{even}}$, of the ring $H^*(X;\Z)/\Tor$:
$$
\alpha_i\alpha_j = c^{k}_{ij}\beta_k,  \ \beta_i\beta_j = d^{k}_{ij}\beta_k,  \ \alpha_i\beta_j = e^{k}_{ij}\alpha_k,  \   c^*_{**}, d^*_{**}, e^*_{**}\in \Z. 
$$
Let us denote by $\bar{\gamma}_1,\bar{\gamma}_2,\ldots$ the elements of the type $(*)$, which forms a $\Z$-basis of our ring $H^*(\Sym^nX;\Z)/\Tor$. Now, we will present a simple algorithm for computing the integral multiplication table for this basis $\bar{\gamma}_i \bar{\gamma}_j = \nu^k_{ij} \bar{\gamma}_k, \nu^*_{**}\in \Z$. We have the following calculation:
$$
\bar{\gamma}_i \bar{\gamma}_j = \frac{1}{r_i!}\frac{1}{r_j!}\left(\sum_{\sigma\in S_n} \sigma^{-1}(\omega_i)\right)\left(\sum_{\tau\in S_n} \tau^{-1}(\omega_j)\right) = \frac{1}{r_i!}\frac{1}{r_j!}\sum_{\sigma,\tau \in S_n} \sigma^{-1}(\omega_i)\tau^{-1}(\omega_j) =
$$
$$
=  \frac{1}{r_i!}\frac{1}{r_j!}\sum_{\sigma,\tau \in S_n} \sigma^{-1}(\omega_i) (\tau\sigma)^{-1}(\omega_j) =  \frac{1}{r_i!}\frac{1}{r_j!}\sum_{\tau \in S_n}\sum_{\sigma \in S_n}\sigma^{-1}(\omega_i\tau^{-1}(\omega_j))   =   \sum_{\tau \in S_n} \mu^k_{\tau,ij}\bar{\gamma}_k  =  \nu^k_{ij}\bar{\gamma}_k.
$$
Here, the constants $\mu^*_{*,**}\in \Q$ are easily computable from the original multiplication table of the ring $H^*(X;\Z)/\Tor$. The constants $\nu^k_{ij}= \sum_{\tau\in S_n}\mu^k_{\tau,ij} \in \Q$ are also computable. We know, that the elements $\bar{\gamma}_1,\bar{\gamma}_2,\ldots$ forms the $\Z$-basis of the ring $H^*(\Sym^nX;\Z)/\Tor$. This fact  implies that these rational structural constants $\nu^{k}_{ij}$ are indeed integral. 

So, we have just presented a simple algorithm for computing the integral multiplication table for the constructed $\Z$-basis $(*)$ of our ring $H^*(\Sym^nX;\Z)/\Tor$. The theorem is completely proved.    \  \  \  \  \  \  \  \  \  \  \  \ \ \  \ \ \  \  \  \  \ $\Box$

\medskip

Let us denote by $\mathcal{R}_Z$ the following category. Objects of $\mathcal{R}_Z$ are connected graded commutative rings with identity $A^* = \Z\left< 1 \right>\oplus A^1\oplus A^2\oplus \ldots $ such that all homogeneous components $A^i, i\ge 1,$ are free abelian groups of finite rank. Morphisms of the category $\mathcal{R}_Z$ are ring homomorphisms. Consider the following functor $S^n_{Z}:\mathcal{R}_Z\to \mathcal{R}_Z, n\ge 2$. 

For any $A^*$ from $\mathcal{R}_Z$ take $S^n A^* := (A^{\otimes n})^{S_n}$ and define $S^n_{Z}A^*$ to be  the minimal subring of $S^n A^*$, that contains elements $\chi(a) = a\otimes 1\otimes \ldots \otimes 1 + \ldots + 1\otimes\ldots\otimes 1\otimes a$ for all $a\in A^*$. For any two $A^*$ and $B^*$ from $\mathcal{R}_Z$ and an arbitrary ring homomorphism $f:A^*\to B^*$ we have ring homomorphisms $f^{\otimes n}: A^{\otimes n} \to B^{\otimes n}$ and $S^n f := f^{\otimes n}|_{S^n A^*, S^n B^*}$. Let us define the ring homomorphism $S^n_Z f := S^n f|_{S^n_{Z} A^*, S^n_{Z} B^*}$. This ring homomorphism is correct due to the fact that $f^{\otimes n}(\chi(a)) = \chi(f(a))$ for all $a\in A^*$. Therefore, we obtain a correct functor $S^n_{Z}:\mathcal{R}_Z\to \mathcal{R}_Z$.

\textbf{Lemma 2.} {\it Fix any ring $A^*$ from $\mathcal{R}_Z$. Let us denote by $\alpha_1,\alpha_2,\ldots$ an arbitrary $\Z$-basis of the odd-dimensional part $A^{\odd}$, and by $\beta_1,\beta_2,\ldots$  an arbitrary $\Z$-basis of the even-dimensional part $A^{\even\ge 2}$. Then the following elements 
$$
m_1!m_2!\ldots m_l!   \left(   \frac{1}{r! m_1! m_2! \ldots m_l!}   \sum_{  \sigma\in S_n  } \sigma^{-1}(\alpha_{i_1}\otimes \alpha_{i_2}\otimes \ldots \otimes  \alpha_{i_k} \otimes \underbrace{\beta_{j_1} \otimes \ldots \otimes \beta_{j_1}}_{ \mbox{$m_1$ \text{times} } } \otimes \right. 
$$
$$
\left. \underbrace{\beta_{j_2} \otimes \ldots \otimes \beta_{j_2}}_{\mbox{$m_2$ \text{times} } } \otimes \ldots  \otimes \underbrace{ \beta_{j_l} \otimes \ldots \otimes \beta_{j_l} }_{ \mbox{$m_l$ \text{times} } }\otimes \underbrace{  1\otimes \ldots \otimes 1 }_{\mbox{$r$ \text{times} }}   ) \right)  \in S^nA^*,  \  \eqno{(**)}
$$
where $k\ge 0, l\ge 0, k+l\ge 1, 1\le  i_1<i_2<\ldots <i_k \le B_{\odd}:= \mathrm{rk}{A^{\odd}}, 1\le j_1<j_2<\ldots <j_l \le \tilde{B}_{\even}:= \mathrm{rk}{A^{\even \ge 2}}, m_1\ge 1,m_2\ge 1,\ldots, m_l\ge 1,  r\ge 0, k+m_1+m_2+\ldots + m_l + r = n$, \\
lie in $S^n_{Z}A^*$ and form a $\Z$-basis of this ring.}

\textbf{Proof.} A standard exercise shows that the elements $(**)$ are $\Q$-linearly independent in $(S^n A^*)\otimes \Q = S^n (A^*\otimes \Q)$. So, to conclude the proof of the lemma we only need to prove that these elements lie in $S^n_{Z}A^*$ and $S^n_{Z}A^* = \Z\left<(**)\right>$. 

Take any homogeneous elements $a_1, a_2, \ldots, a_k\in A^{*\ge 1}, 1\le k\le n,$ and consider the following tensor
$$
\chi(a_1,\ldots,a_k) :=  \frac{1}{(n-k)!} \sum_{\sigma\in S_n} \sigma^{-1} (a_1\otimes a_2\otimes\ldots \otimes a_k\otimes 1\otimes \ldots \otimes 1) \in S^n A^*.   \eqno{(***)}
$$
By the same argumentation as in the proof of Theorem 1, one can show that $\chi(a_1,\ldots,a_k)\in S^n_{Z}A^*$. It is also easy to see that $\Z\left<(**)\right> = \Z\left<(***)\right>$. Therefore, we have the double inclusion 
$$
\Z\left< \chi(a), a\in A^{*\ge 1} \right> \subset \Z\left<(***)\right> \subset S^n_{Z}A^* = \Z\left< \chi(b_1)\chi(b_2)\ldots \chi(b_N), N\ge 1, b_1,\ldots,b_N\in A^{*\ge 1}  \right>. 
$$ 
Now, to prove the desired coincidence $\Z\left<(***)\right> = S^n_{Z}A^*$ we only need to check that 
$$
\chi(a_1,\ldots,a_k)\chi(b)\in \Z\left<(***)\right>
$$
for any $a_1,\ldots,a_k, b \in A^{*\ge 1}, 1\le k\le n$. Here is the calculation: 
$$
\chi(a_1,\ldots,a_k)\chi(b) = \frac{1}{(n-k)!} \frac{1}{(n-1)!} \sum_{\sigma, \tau \in S_n} \sigma^{-1} (a_1\otimes a_2\otimes\ldots \otimes a_k\otimes 1\otimes \ldots \otimes 1) \tau^{-1}(b\otimes 1\otimes \ldots \otimes 1) =
$$
$$
= \frac{1}{(n-k)!} \frac{1}{(n-1)!} \sum_{\sigma, \tau \in S_n} \sigma^{-1} (a_1\otimes a_2\otimes\ldots \otimes a_k\otimes 1\otimes \ldots \otimes 1) \sigma^{-1} \tau^{-1}(b\otimes 1\otimes \ldots \otimes 1) = 
$$
$$
=  \frac{1}{(n-k)!} \frac{1}{(n-1)!} \sum_{\sigma, \tau \in S_n} \sigma^{-1} ((a_1\otimes a_2\otimes\ldots \otimes a_k\otimes 1\otimes \ldots \otimes 1) \tau^{-1}(b\otimes 1\otimes \ldots \otimes 1)) = 
$$
$$
= \  \frac{1}{(n-k)!} \sum_{\sigma \in S_n}  \sigma^{-1} ( a_1\otimes a_2\otimes\ldots \otimes a_k\otimes 1\otimes \ldots \otimes 1 \cdot  b\otimes 1\otimes \ldots \otimes 1 ) \ +
$$
$$
+ \     \frac{1}{(n-k)!} \sum_{\sigma \in S_n}  \sigma^{-1} ( a_1\otimes a_2\otimes\ldots \otimes a_k\otimes 1\otimes \ldots \otimes 1 \cdot  1\otimes b\otimes 1\otimes \ldots \otimes 1 )  \ +   \  \ldots  \  + 
$$
$$
+ \   \frac{1}{(n-k)!} \sum_{\sigma \in S_n}  \sigma^{-1} ( a_1\otimes a_2\otimes\ldots \otimes a_k\otimes 1\otimes \ldots \otimes 1 \cdot  1\otimes 1\otimes \ldots \otimes b ) \ = 
$$
$$
= \ (-1)^{|b|(|a_2|+\ldots +|a_k|)}\chi(a_1b,a_2,\ldots,a_k) \ + \  (-1)^{|b|(|a_3|+\ldots +|a_k|)}\chi(a_1,a_2b,a_3,\ldots,a_k)  \  +  \  \ldots    \  + 
$$
$$
+   \   \chi(a_1,a_2,\ldots,a_kb)  \  +   \   \chi(a_1,a_2,\ldots,a_k, b).
$$
(In this formula if $k=n$, then by definition $\chi(a_1,a_2,\ldots,a_k, b)=0$).    \ \ \ \ \ \ \ \ \  \ \ \ \  \ \ \ \ \ \ \ \ \ \ \ \   \ \ \ \ \ \ \ \ \  \ \ \ \  \ \ \ \ \ \ \ \ \ \ \ \  \  \  \  \ \ \ \    $\Box$

\medskip 

Lemma 2 easily implies that \\
(i) $S^n_{Z}A^* = S^n A^*$ {\it iff} $A^{2m}=0, \forall m\ge 1$; \\
(ii) $S^n_{Z}A^i = S^n A^i$ for $i=1,2,3$; \\
(iii) if $A^{2m}\ne 0$ for some $m\ge 1$, then for any $a\in A^{2m}, a\ne s a', s\ge 2,$ the tensor $Na\otimes a\otimes\ldots \otimes a$ belongs to $S^n_{Z}A^*$ {\it iff} $n! | N$; \\
(iv) if $A^{2m}\ne 0$ for some $m\ge 1$, then $n!S^n A^* \subset S^n_{Z}A^*$. 

Evidently, Theorem 1 implies 

\textbf{Corollary 2.}  {\it Let $X$ be a connected Hausdorff space homotopy equivalent to a CW-complex and having finitely generated integral homology groups in all dimensions. Fix any $n\ge 2$. Denote by $A^*$ the ring $H^*(X;\Z)/\Tor$. Then the ring $H^*(\Sym^n X;\Z)/\Tor$ is just the functor $S^n_{Z} A^*$. For any two spaces $X$ and $Y$ satisfying the above mentioned regularity conditions and any continuous map $f\colon Y \to X$ the corresponding map $(\Sym^n f)^*_{/\Tor} \colon H^*(\Sym^n X;\Z)/\Tor \to H^*(\Sym^n Y;\Z)/\Tor$ is just $S^n_{Z} (f^*_{/\Tor})$, where $f^*_{/\Tor} \colon H^*(X;\Z)/\Tor \to H^*(Y;\Z)/\Tor$ is the induced homomorphism.}

\section{Macdonald's theorem}

The natural homeomorphism $\Sym^n \Cn\cong \Cn^n$ implies that for any 2-dimensional manifold $M^2$ its symmetric product $\Sym^nM^2$ is a $2n$-dimensional topological manifold. Actually, it is always smoothable and admits a natural complex structure if $M^2$ is equipped with a structure of a Riemann surface. 

It is well known that if $M^2$ is a compact Riemann surface then $\Sym^n M^2$ is a smooth projective algebraic variety. The standard fact is that $\Sym^n \Cn P^1 \cong \Cn P^n$. Let us denote by $M^2_g$ an arbitrary compact Riemann surface of genus $g\ge 1$. There is a classical Abel map $Ab:\Sym^n M^2_g \to Jac=T^g$, where $Jac$ is a Jacobian variety of $M^2_g$, which is a $g$-dimensional complex torus. The following famous theorem states the good behavior of Abel map for stable $n$.

{\bf Theorem\,(Mattuck,\,1961).} {\it If $n\ge 2g-1$, then Abel map $Ab:\Sym^n M^2_g \to Jac=T^g$ is an algebraic projective bundle, with fiber $\mathbb{C}P^{n-g}$. Moreover, it is a projectivization of a canonically defined algebraic complex vector bundle, with fiber $\mathbb{C}^{n-g+1}$.}

The integral cohomology ring of $M^2_g$ is well known:
$$
H^*(M^2_g;\Z)=\Z \left< 1\right>\oplus \Z\left< \alpha_1,\ldots,\alpha_{2g}\right>\oplus \Z \left< \beta \right>,
$$
$$
\alpha_i\alpha_{i+g}= - \alpha_{i+g}\alpha_i=\beta, 1\le i \le g, \text{the rest products are zero.} 
$$

Let us denote
$$
\chi(\alpha_i)=\alpha_i\otimes 1\otimes\ldots\otimes 1 +\ldots + 1\otimes\ldots\otimes 1\otimes\alpha_i= \xi_i, 1\le i\le g,
$$
$$
\chi(\alpha_{i+g})=\alpha_{i+g}\otimes 1\otimes\ldots\otimes 1 +\ldots + 1\otimes\ldots\otimes 1\otimes\alpha_{i+g}= \xi'_i, 1\le i\le g,
$$
$$
\chi(\beta)=\beta\otimes 1\otimes\ldots\otimes 1 +\ldots + 1\otimes\ldots\otimes 1\otimes\beta= \eta.
$$

Without Integrality Lemma we have
$$
\xi_1,\ldots,\xi_g,\xi'_1,\ldots,\xi'_g,\eta \in H^*(\Sym^n M^2_g;\Q)
$$

Consider the free graded commutative $\Q$-algebra
$$
\Lambda_{\Q} \left< x_1,\ldots,x_g,x'_1,\ldots,x'_g \right>\otimes \Q[y],
$$
where $x_1,\ldots,x_g,x'_1,\ldots,x'_g,y$ are formal variables, $|x_i|=|x'_i|=1, 1\le i\le g, |y|=2$. 

Here we state a famous theorem of Macdonald \cite{Mac1} (case $\Q$). 

{\bf Theorem\,$\alpha$\,(Macdonald,\,1962).} {\it Suppose $M^2_g$ is an arbitrary compact Riemann surface of genus $g\ge 1$, and $n\ge 2$. Take the $\Q$-algebra homomorphism
$$
f_{\Q}: \Lambda_{\Q} \left< x_1,\ldots,x_g,x'_1,\ldots,x'_g \right>\otimes \Q[y] \to H^*(\Sym^n M^2_g;\Q),
$$
$$
x_i \mapsto \xi_i, \ \  x'_i \mapsto \xi'_i, \ \ 1\le i\le g, \ \ y \mapsto \eta.
$$
Then the following statements hold: 

(i) $f_{\Q}$ is an epimorphism, and $\mathrm{Ker}(f_{\Q}) = I^*_{Mac}$ is generated by the following polynomials
$$
x_{i_1}\ldots x_{i_a}x'_{j_1}\ldots x'_{j_b} (y - x_{k_1}x'_{k_1})\ldots  (y - x_{k_c}x'_{k_c}) y^q,
$$
where $a+b+2c+ q=n+1$ and $i_1,\ldots,i_a,j_1,\ldots,j_b,k_1,\ldots,k_c$ are distinct integers from $1$ to $g$ inclusive.

(ii) if $n\ge 2g-1$, then $I^*_{Mac}$ is generated by the single polynomial
$$
(y - x_1x'_1)(y - x_2x'_2)\ldots (y - x_g x'_g) y^{n-2g+1}.
$$
It is the only polynomial with $q=q_{min}=n-2g+1$.

(iii) if $2\le n\le 2g-2$, then $I^*_{Mac}$ is generated by $\binom{2g}{n+1}$ polynomials with $q=0$.}

The statements $(i)$ and $(ii)$ hold true, but the statement $(iii)$ needs the following correction:\\
$(\widetilde{iii}_A)$ if $n=3,5,\dots,2g-3$, then $I^*_{Mac}$ is generated by $\binom{2g}{n+1}$ polynomials with $q=0$.\\
$(\widetilde{iii}_B)$ if $n=2,4,\dots,2g-2$, then $I^*_{Mac}$ is generated by $\binom{2g}{n+1}$ polynomials with $q=0$ and one more polynomial $(y-x_1x'_1)(y-x_2x'_2)\ldots(y-x_{\frac{n}{2}}x'_{\frac{n}{2}})y$.

In Macdonald's paper the proofs of the statements $(i)$ and $(ii)$ are correct. Here we reprove the statement $(ii)$ and prove the correction $(\widetilde{iii}_A)$ and $(\widetilde{iii}_B)$. We will deduce these statements from $(i)$. 

In fact a more weaker correction of the statement $(iii)$ was already made by Bertram and Thaddeus in 2001 (see \cite{Bert}, the remark after Theorem 2.2). Namely, they just mentioned that the original proof by Macdonald gives that the statement $(iii)$ could be changed by the following true statement: 

\noindent $(\widetilde{iii}_C)$ if $2\le n\le 2g-2$, then $I^*_{Mac}$ is generated by $\binom{2g}{n+1}$ polynomials with $q=0$ and $\binom{2g}{n}$ polynomials with $q=1$.

\noindent This was pointed out to me by one of the referees, and I am very grateful to him.

{\bf Theorem 2.} {\it Suppose $M^2_g$ is an arbitrary compact Riemann surface of genus $g\ge 1$, and $n\ge 2$. Take the $\Q$-algebra epimorphism
$$
f_{\Q}: \Lambda_{\Q} \left< x_1,\ldots,x_g,x'_1,\ldots,x'_g \right>\otimes \Q[y] \to H^*(\Sym^n M^2_g;\Q),
$$
$$
x_i \mapsto \xi_i, \ \  x'_i \mapsto \xi'_i, \ \ 1\le i\le g, \ \ y \mapsto \eta.
$$
Let us denote by $I^*_{Mac}$ the ideal $\mathrm{Ker}(f_{\Q})$. Then the following statements hold true:\\
$(\widetilde{iii}_A)$ if $n=3,5,\dots,2g-3$, then $I^*_{Mac}$ is generated by $\binom{2g}{n+1}$ polynomials with $q=0$. None of these polynomials can be removed.\\
$(\widetilde{iii}_B)$ if $n=2,4,\dots,2g-2$, then $I^*_{Mac}$ is generated by $\binom{2g}{n+1}$ polynomials with $q=0$ and one more polynomial $(y-x_1x'_1)(y-x_2x'_2)\ldots(y-x_{\frac{n}{2}}x'_{\frac{n}{2}})y$.  None of these polynomials can be removed.}

{\bf Proof.} We will deduce these statements from $(i)$. It is evident that the 4-tuples $(a,b,c,q)\in \Z_{+}^4$ from $(i)$ are the solutions of the following system
$$
\left\{ \begin{aligned}
a+b+2c+q & =n+1 \\
a+b+c & \le  g .
\end{aligned}
\right.
$$

Now, our goal is to deduce any polynomial  from $(i)$ with the 4-tuple $(\tilde{a},\tilde{b},\tilde{c},\tilde{q})$ from polynomials with $(a,b,c,q)$ such that $q<\tilde{q}$ (if it can be done). Consider the following 4 cases.

{\bf Case I. $(\tilde{q}\ge 1, \tilde{a}\ge 1)$}\\
Set $a=\tilde{a}-1, b=\tilde{b}, c=\tilde{c}+1, q=\tilde{q}-1$. Take any distinct integers 
$$
\{ i_1,i_2,\ldots,i_a,i_{\tilde{a}=a+1},j_1,\ldots,j_{\tilde{b}=b}, k_1,\ldots,k_{\tilde{c}=c-1} \}\subset \{1,2,\ldots,g\}.
$$
We will deduce the polynomial 
$$
x_{i_1}\ldots x_{i_{\tilde{a}}}x'_{j_1}\ldots x'_{j_{\tilde{b}}} (y - x_{k_1}x'_{k_1})\ldots  (y - x_{k_{\tilde{c}}}x'_{k_{\tilde{c}}}) y^{\tilde{q}} \in I^*_{Mac}
$$
from the polynomial for the set of indices 
$$
\{ i_1,i_2,\ldots,i_a,j_1,\ldots,j_{b}, k_1,\ldots,k_{c-1}, k_{c}:=i_{a+1} \}\subset \{1,2,\ldots,g\},
$$
and $q=\tilde{q}-1$. Here is the deduction:
$$
x_{i_1}\ldots x_{i_a}x'_{j_1}\ldots x'_{j_b} (y - x_{k_1}x'_{k_1})\ldots  (y - x_{k_c}x'_{k_c}) y^q \in I^*_{Mac}  \  \Rightarrow
$$
$$
x_{i_1}\ldots x_{i_{\tilde{a}-1}}x'_{j_1}\ldots x'_{j_{\tilde{b}}} (y - x_{k_1}x'_{k_1})\ldots  (y - x_{k_{\tilde{c}}}x'_{k_{\tilde{c}}})(y - x_{i_{\tilde{a}}}x'_{i_{\tilde{a}}}) y^{\tilde{q}-1} \in I^*_{Mac} \ \Rightarrow
$$
$$
x_{i_1}\ldots x_{i_{\tilde{a}-1}}x'_{j_1}\ldots x'_{j_{\tilde{b}}} (y - x_{k_1}x'_{k_1})\ldots  (y - x_{k_{\tilde{c}}}x'_{k_{\tilde{c}}})y^{\tilde{q}}     =  
$$
$$
=    x_{i_1}\ldots x_{i_{\tilde{a}-1}} x_{i_{\tilde{a}}}x'_{i_{\tilde{a}}} x'_{j_1}\ldots x'_{j_{\tilde{b}}} (y - x_{k_1}x'_{k_1})\ldots  (y - x_{k_{\tilde{c}}}x'_{k_{\tilde{c}}})y^{\tilde{q}-1} \  \text{mod} \ I^*_{Mac} \ \Rightarrow  \  (\mbox{multiply both parts by} \  x_{i_{\tilde{a}}})  
$$
$$
\Rightarrow \  x_{i_1}\ldots x_{i_{\tilde{a}}}x'_{j_1}\ldots x'_{j_{\tilde{b}}} (y - x_{k_1}x'_{k_1})\ldots  (y - x_{k_{\tilde{c}}}x'_{k_{\tilde{c}}}) y^{\tilde{q}} \in I^*_{Mac}.
$$

{\bf Case II. $(\tilde{q}\ge 1, \tilde{b}\ge 1)$}\\
Set $a=\tilde{a}, b=\tilde{b}-1, c=\tilde{c}+1, q=\tilde{q}-1$. The reasoning in this case is absolutely similar and is left to the reader. 

{\bf Case III. $(\tilde{q}\ge 2, \tilde{a} = \tilde{b} = 0, 0\le \tilde{c}\le g-1)$}\\
Set $a=b=0, c=\tilde{c}+1, q=\tilde{q}-2$. Take any distinct integers 
$$
\{ k_1,\ldots,k_{\tilde{c}=c-1} \}\subset \{1,2,\ldots,g\},
$$
and add to them any $k_c\in \{1,2,\ldots,g\} \backslash   \{ k_1,\ldots,k_{\tilde{c}=c-1} \}$. Here is the desired deduction:
$$
(y - x_{k_1}x'_{k_1})\ldots  (y - x_{k_c}x'_{k_c}) y^q \in I^*_{Mac} \  \Rightarrow
$$
$$
(y - x_{k_1}x'_{k_1})\ldots  (y - x_{k_{c-1}}x'_{k_{c-1}}) y^q x_{k_c}x'_{k_c}= (y - x_{k_1}x'_{k_1})\ldots  (y - x_{k_{c-1}}x'_{k_{c-1}}) y^{q+1}  \ \text{mod} \ I^*_{Mac} \   \Rightarrow
$$
$$
[A:=(y - x_{k_1}x'_{k_1})\ldots  (y - x_{k_{c-1}}x'_{k_{c-1}}) y^q]  \  \Rightarrow \  Ax_{k_c}x'_{k_c} = Ay \ \text{mod} \ I^*_{Mac} \  \Rightarrow 
$$
$$
I^*_{Mac}\ni Ax_{k_c}x'_{k_c}x_{k_c}x'_{k_c} = Ax_{k_c}x'_{k_c}y=Ay^2 \  \Rightarrow  \  (y - x_{k_1}x'_{k_1})\ldots  (y - x_{k_{\tilde{c}}}x'_{k_{\tilde{c}}}) y^{\tilde{q}} \in I^*_{Mac}.
$$

{\bf Case IV. $(\tilde{a} = \tilde{b} = 0, \tilde{c} = g)$}\\ 
As we have $\tilde{a} + \tilde{b} + 2\tilde{c} + \tilde{q} = 2g + \tilde{q} = n+1$ and $\tilde{q}\ge 0$, this case can occur {\it iff} $n\ge 2g-1$ (stable dimension). In this case there is only one polynomial:
$$
(y - x_1x'_1)(y - x_2x'_2)\ldots (y - x_g x'_g) y^{n-2g+1} \in I^*_{Mac}.
$$

Suppose first that $2\le n\le 2g-2$. Then from above analysis of Cases I,II,III, we obtain the following fact: the ideal $I^*_{Mac}$ is generated by the polynomials from $(i)$ with $q=0$ and $q=1$. 

Consider any 4-tuple $(\tilde{a},\tilde{b},\tilde{c},\tilde{q}=1)$. We have $\tilde{a} + \tilde{b} + 2\tilde{c} + 1 = n+1$. If $n$ is an odd number, $n=3,5,\ldots,2g-3$, then $\tilde{a} + \tilde{b}$ is always odd. In particular, $\tilde{a}\ge 1$ or $\tilde{b}\ge 1$. So, we come to Case I or II. Thus, any polynomial with a given $(\tilde{a},\tilde{b},\tilde{c},\tilde{q}=1)$ is a consequence of polynomials with $q=0$. The correction $(\widetilde{iii}_A)$ is completely proved.

Now, suppose that $n=2,4,\ldots,2g-2$. It is easy to check, that in this case the ideal $I^*_{Mac}$ is generated by the polynomials from $(i)$ with $q=0$ and $(a,b,c,q)=(0,0,\frac{n}{2},1)$. Any polynomial with 4-tuple $(0,0,\frac{n}{2},1)$ has the form:
$$
P_{k_1,k_2,\ldots,k_{\frac{n}{2}}} := y (y-x_{k_1}x'_{k_1}) (y-x_{k_2}x'_{k_2})\ldots (y-x_{k_{\frac{n}{2}}}x'_{k_{\frac{n}{2}}}) \in I^*_{Mac},
$$
where $\{ k_1,k_2,\ldots,k_{\frac{n}{2}} \} \subset \{1,2,\ldots,g\} $ and $|\{ k_1,k_2,\ldots,k_{\frac{n}{2}} \}| = \frac{n}{2}$.

Denote by $I^*_{\left< q=0 \right>}$ the ideal of our ring which is generated by all polynomials with $q=0$. Then for any distinct integers 
$$
i,k,k_2,k_3,\ldots,k_{\frac{n}{2}}\in \{1,2,\ldots,g \}
$$
we have
$$
x_i(y-x_kx'_k)(y-x_{k_2}x'_{k_2})\ldots (y-x_{k_{\frac{n}{2}}}x'_{k_{\frac{n}{2}}}) \in  I^*_{\left< q=0 \right>}.
$$
Multiplying by $x'_i$ we get
$$
x_ix'_i(y-x_kx'_k)(y-x_{k_2}x'_{k_2})\ldots (y-x_{k_{\frac{n}{2}}}x'_{k_{\frac{n}{2}}}) \in  I^*_{\left< q=0 \right>},
$$
$$
y x_ix'_i(y-x_{k_2}x'_{k_2})\ldots (y-x_{k_{\frac{n}{2}}}x'_{k_{\frac{n}{2}}}) = x_kx'_k x_ix'_i (y-x_{k_2}x'_{k_2})\ldots (y-x_{k_{\frac{n}{2}}}x'_{k_{\frac{n}{2}}}) \ \mbox{mod} \  I^*_{\left< q=0 \right>}.
$$ 
By change $i\leftrightarrow k$ we obtain
$$
y x_kx'_k(y-x_{k_2}x'_{k_2})\ldots (y-x_{k_{\frac{n}{2}}}x'_{k_{\frac{n}{2}}}) = x_ix'_i x_kx'_k (y-x_{k_2}x'_{k_2})\ldots (y-x_{k_{\frac{n}{2}}}x'_{k_{\frac{n}{2}}}) \ \mbox{mod} \  I^*_{\left< q=0 \right>}.
$$ 
So, for any distinct integers $i,k,k_2,k_3,\ldots,k_{\frac{n}{2}}\in \{1,2,\ldots,g \}$ we have
$$
y (x_ix'_i - x_kx'_k)(y-x_{k_2}x'_{k_2})\ldots (y-x_{k_{\frac{n}{2}}}x'_{k_{\frac{n}{2}}}) = 0  \ \mbox{mod} \  I^*_{\left< q=0 \right>}.
$$
This formula can be rewritten in the way
$$
P_{i,k_2,\ldots,k_{\frac{n}{2}}} =  P_{k,k_2,\ldots,k_{\frac{n}{2}}} \ \mbox{mod} \  I^*_{\left< q=0 \right>}.
$$
The last formula implies that for any $1\le k_1< k_2<\ldots<k_{\frac{n}{2}}\le g$ we get the identity
$$
P_{k_1,k_2,k_3,\ldots,k_{\frac{n}{2}}} = P_{1,k_2,k_3,\ldots,k_{\frac{n}{2}}} = P_{1,2,k_3,\ldots,k_{\frac{n}{2}}} = \ldots = P_{1,2,\ldots,\frac{n}{2}}  \ \mbox{mod} \  I^*_{\left< q=0 \right>}.
$$

Therefore, in the case $n=2,4,\ldots,2g-2$, the ideal $I^*_{Mac}$ is generated by the $\binom{2g}{n+1}$ polynomials with $q=0$ and one more polynomial $P_{1,2,\ldots,\frac{n}{2}} = (y-x_1x'_1)(y-x_2x'_2)\ldots(y-x_{\frac{n}{2}}x'_{\frac{n}{2}})y$ with $q=1$. Moreover, the polynomial $P_{1,2,\ldots,\frac{n}{2}}$ does not belong to the ideal $I^*_{\left< q=0 \right>}$. It can be proved by the following argumentation.

We have $\mathrm{deg}P_{1,2,\ldots,\frac{n}{2}} = n+2,$ and $\mathrm{deg}R = n+1$ for any our polynomial $R$ with $q=0$. So, if $P_{1,2,\ldots,\frac{n}{2}} \in I^*_{\left< q=0 \right>}$, then there should be the identity
$$
P_{1,2,\ldots,\frac{n}{2}} = y^{\frac{n}{2} + 1} + \ldots = \sum_{i=1}^g x_i R_i + \sum_{j=1}^g x'_jR'_j,
$$
for some homogeneous elements $R_1,\ldots,R_g,R'_1,\ldots,R'_g \in \Lambda_{\Q} \left< x_1,\ldots,x_g,x'_1,\ldots,x'_g \right>\otimes \Q[y]$ of degree $n+1$. But, as the ring $\Lambda_{\Q} \left< x_1,\ldots,x_g,x'_1,\ldots,x'_g \right>\otimes \Q[y]$ is a free graded commutative ring, the above identity cannot hold true. The correction $(\widetilde{iii}_B)$ is completely proved.

Let us consider the stable case $n\ge 2g -1.$ Take any 4-tuple $(a,b,c,q)$. We have double inequality $a+b+2c\le 2(a+b+c)\le 2g$. We also have the identity $a+b+2c+q=n+1$. It follows that $q= n+1 - (a+b+2c)\ge n+1 - 2g=q_{min}\ge 0$. Thus, the minimal $q$ that may occur in this case is $q_{min}=n-2g+1$ and there is only one polynomial with this $q_{min}$:
$$
(y-x_1x'_1)(y-x_2x'_2)\ldots (y-x_gx'_g)y^{n-2g+1} \in I^*_{Mac}.
$$ 

Consider any 4-tuple $(\tilde{a},\tilde{b},\tilde{c},\tilde{q})$ with $\tilde{q}>q_{min}$. If $\tilde{a}\ge 1$ or $\tilde{b}\ge 1$, then we get to Cases I or II, and, therefore, any polynomial with this 4-tuple $(\tilde{a},\tilde{b},\tilde{c},\tilde{q})$ is a consequence of some polynomial with $q<\tilde{q}$. If $\tilde{a}=\tilde{b}=0$, then $\tilde{q}=n+1-2\tilde{c}> n+1 - 2g$. It follows that $0\le \tilde{c}\le g-1$ and $\tilde{q}\ge 2$, so, we get to Case III. And again any polynomial with this 4-tuple $(\tilde{a},\tilde{b},\tilde{c},\tilde{q})$ is a consequence of some polynomial  with $q<\tilde{q}$. 

By finite descending induction we obtain, that in the stable case $n\ge 2g -1$ the ideal $I^*_{Mac}$ is generated by the single polynomial with $q=q_{min}$. The part $(ii)$ is completely proved.  

Suppose $2\le n\le 2g-2$. So, we get to the cases $(\widetilde{iii}_A)$ and $(\widetilde{iii}_B)$. We will prove now, that no one of the $\binom{2g}{n+1}$ polynomials with $q=0$ can be removed.

Fix any $2\le n \le 2g-2$. Denote by $B_k, 0\le k\le 2n,$ the $k$-th rational Betti number of $\Sym^n M^2_g$. In Macdonald's paper there is an explicit formula for $B_k$:
$$
B_{2n-k} = B_k = \binom{2g}{k} + \binom{2g}{k-2} + \binom{2g}{k-4} + \ldots, \ \forall 0\le k\le n.
$$
So, we get 
$$
\dim H^{n+1}(\Sym^n M^2_g;\Q) = B_{n+1} = B_{n-1} = \binom{2g}{n-1} + \binom{2g}{n-3} + \ldots.
$$
We also have
$$
(\Lambda_{\Q} \left< x_1,\ldots,x_g,x'_1,\ldots,x'_g \right>\otimes \Q[y])^{n+1}/ I^{n+1}_{Mac} \cong H^{n+1}(\Sym^n M^2_g;\Q).
$$
It is easy to compute
$$
\dim (\Lambda_{\Q} \left< x_1,\ldots,x_g,x'_1,\ldots,x'_g \right>\otimes \Q[y])^{n+1} = \binom{2g}{n+1} + \binom{2g}{n-1} + \binom{2g}{n-3} + \ldots.
$$
It follows that $\dim I^{n+1}_{Mac} = \binom{2g}{n+1}$. But, any polynomial from $(i)$ has degree $a+b+2c+2q = n + 1 + q\ge n+1$. So, $\Q$-vector space $I^{n+1}_{Mac}$ of dimension $\binom{2g}{n+1}$ is linearly generated by $\binom{2g}{n+1}$ polynomials with $q=0$. Therefore, the system of $\binom{2g}{n+1}$ polynomials with $q=0$ is $\Q$-linearly independent, so, none of them can be removed. The theorem is completely proved.  \ \ \ \ \ \  \ \ \ \ \ \ \ \ \ \ \ \ \ \ \ \ \ \ \ \ \ \ \ \ \ \ \ \ \ \ \ \ \ \ \ \ \ \ \ \ \ \ \ \ \ \ \ \ \ \ \ \ \ \ \ \ \ \ \ \  \ \ \ \ \ \ \ \ \ \  \ \ \ \ \ \ \ \ \ \  \ \ \ \ \ \ \ \ \ \ \ \ \ \ \ \ \   $\Box$

\medskip

In the same paper \cite{Mac1} Macdonald proves that the integral cohomology ring $H^*(\Sym^n M^2_g;\Z)$ has no torsion. This fact can be derived more directly in the following way. Suppose $X$ is a connected CW-complex. Let us fix any base-point $x_0\in X$. This gives the series of inclusion mappings 
$$
X\hookrightarrow \Sym^2 X \hookrightarrow \Sym^3 X \hookrightarrow \ldots.
$$ 
The prominent Steenrod Splitting Theorem states the isomorphism of abelian groups 
$$
H_q(\Sym^n X;\Z)\cong \bigoplus_{k=1}^{n} H_q(\Sym^k X/\Sym^{k-1} X;\Z), \ q\ge 1, \ n\ge 1,
$$ 
and also that $H_q(\Sym^n X;\Z)$ is mapped isomorphically to a direct summand of $H_q(\Sym^{n+1} X;\Z),  q, n\ge 1$.

\noindent This theorem implies that if $X$ is a connected finite CW-complex and integral homology $H_*(\Sym^N X;\Z)$ is torsion-free for all $N$ greater than some $N_0$, then integral homology $H_*(\Sym^n X;\Z)$ is torsion-free for all $n\ge 1$. 
\noindent If we take $X=M^2_g$, then for $N\ge N_0:= 2g-1$ the homology $H_*(\Sym^N M^2_g;\Z)$ is torsion-free due to the above Mattuck's theorem. Therefore, we get the result that the integral cohomology ring $H^*(\Sym^n M^2_g;\Z)$ has no torsion.

Let us take the following free graded commutative ring
$$
\Lambda_{\Z} \left< x_1,\ldots,x_g,x'_1,\ldots,x'_g \right>\otimes \Z[y],
$$
where $x_1,\ldots,x_g,x'_1,\ldots,x'_g,y$ are formal variables, $|x_i|=|x'_i|=1, 1\le i\le g, |y|=2$. 

Here we formulate the famous theorem of Macdonald \cite{Mac1} (case $\Z$). 

{\bf Theorem\,$\beta$\,(Macdonald,\,1962).} {\it Suppose $M^2_g$ is an arbitrary compact Riemann surface of genus $g\ge 1$, and $n\ge 2$. Take the graded ring homomorphism
$$
f_{\Z}: \Lambda_{\Z} \left< x_1,\ldots,x_g,x'_1,\ldots,x'_g \right>\otimes \Z[y] \to H^*(\Sym^n M^2_g;\Z),
$$
$$
x_i \mapsto \xi_i, \ \  x'_i \mapsto \xi'_i, \ \ 1\le i\le g, \ \ y \mapsto \eta.
$$
Then the following statements hold: 

(i) $f_{\Z}$ is an epimorphism, and $\mathrm{Ker}(f_{\Z}) = I^*_{Mac}$ is generated by the following polynomials
$$
x_{i_1}\ldots x_{i_a}x'_{j_1}\ldots x'_{j_b} (y - x_{k_1}x'_{k_1})\ldots  (y - x_{k_c}x'_{k_c}) y^q,
$$
where $a+b+2c+ q=n+1$ and $i_1,\ldots,i_a,j_1,\ldots,j_b,k_1,\ldots,k_c$ are distinct integers from $1$ to $g$ inclusive.

(ii) if $n\ge 2g-1$, then $I^*_{Mac}$ is generated by the single polynomial
$$
(y - x_1x'_1)(y - x_2x'_2)\ldots (y - x_g x'_g) y^{n-2g+1}.
$$
It is the only polynomial with $q=q_{min}=n-2g+1$.

(iii) if $2\le n\le 2g-2$, then $I^*_{Mac}$ is generated by $\binom{2g}{n+1}$ polynomials with $q=0$.}

The statements $(i)$ and $(ii)$ again hold true, but the statement $(iii)$ needs the same correction, as in the case $\Q$:\\
$(\widetilde{iii}_A)$ if $n=3,5,\dots,2g-3$, then $I^*_{Mac}$ is generated by $\binom{2g}{n+1}$ polynomials with $q=0$.\\
$(\widetilde{iii}_B)$ if $n=2,4,\dots,2g-2$, then $I^*_{Mac}$ is generated by $\binom{2g}{n+1}$ polynomials with $q=0$ and one more polynomial $(y-x_1x'_1)(y-x_2x'_2)\ldots(y-x_{\frac{n}{2}}x'_{\frac{n}{2}})y$.

The deduction of parts $(ii), (\widetilde{iii}_A)$ and $(\widetilde{iii}_B)$ from part $(i)$ is absolutely the same as in the above case $\Q$ (when we derived above the relations of $(i)$ from the corresponding relations of $(ii), (\widetilde{iii}_A)$ and $(\widetilde{iii}_B)$, there were no denominators: all constants were equal to $0$ and $\pm1$).

{\bf The original proof of Theorem $\beta (i)$ was the following claim:}\\
``We have proved Theorem $\alpha (i)$. The integral cohomology ring $H^*(\Sym^n M^2_g;\Z)$ has no torsion. So, everything remains true, when we replace $\Q$ by $\Z$ throughout.''

But, as we have shown in the Introduction, there exist manifolds $L, M$ with torsion-free cohomology, having equal rational cohomology rings $H^*(L;\Q)\cong H^*(M;\Q)$ but nonisomorphic integral cohomology rings $H^*(L;\Z)\ncong H^*(M;\Z)$. 

More careful analysis for the deduction of the result over $\Z$ from the result over $\Q$ gives us the following three gaps. 

{\bf Gap 1.} Why the considered elements
$$
\xi_1,\ldots,\xi_g,\xi'_1,\ldots,\xi'_g,\eta \in H^*(\Sym^n M^2_g;\Q)
$$
lie in the integral lattice $H^*(\Sym^n M^2_g;\Z)$? This Gap is filled in by Integrality Lemma. 

{\bf Gap 2.} Why the graded ring homomorphism 
$$
f_{\Z}: \Lambda_{\Z} \left< x_1,\ldots,x_g,x'_1,\ldots,x'_g \right>\otimes \Z[y] \to H^*(\Sym^n M^2_g;\Z),
$$
$$
x_i \mapsto \xi_i, \ \  x'_i \mapsto \xi'_i, \ \ 1\le i\le g, \ \ y \mapsto \eta,
$$
is an epimorphism? Equivalently, why the elements $\xi_1,\ldots,\xi_g,\xi'_1,\ldots,\xi'_g,\eta$ are multiplicative generators of the ring $H^*(\Sym^n M^2_g;\Z)$? This Gap was filled in by Seroul in \cite{Ser}. It also can be filled in by Theorem 1.

Now, we have the epimorphism
$$
f_{\Z}: \Lambda_{\Z} \left< x_1,\ldots,x_g,x'_1,\ldots,x'_g \right>\otimes \Z[y] \to H^*(\Sym^n M^2_g;\Z),
$$
whose kernel contains $I^*_{Mac}$ (by $I^*_{Mac}$ we denote the ideal generated by the polynomials from $(i)$). Therefore, the induced graded ring homomorphism
$$
\bar{f}_{\Z}: \Lambda_{\Z} \left< x_1,\ldots,x_g,x'_1,\ldots,x'_g \right>\otimes \Z[y]/ I^*_{Mac} \to H^*(\Sym^n M^2_g;\Z)
$$
is also an epimorphism. Moreover, after $\otimes \Q$ it turns out to be an isomorphism (Theorem $\alpha (i)$). It follows that 
$$
\mathrm{Ker}(\bar{f}_{\Z}) = \Tor(\Lambda_{\Z} \left< x_1,\ldots,x_g,x'_1,\ldots,x'_g \right>\otimes \Z[y]/ I^*_{Mac}).
$$

{\bf Gap 3.} Why $\mathrm{Ker}(\bar{f}_{\Z})$ is zero? This Gap was filled in by Seroul in \cite{Ser} (see Theorem 1.2.3, the proof takes 7 pages). Here we present a shorter 3 pages  proof of the statement.

{\bf Proposition 1.} {\it In the above notations the kernel $\mathrm{Ker}(\bar{f}_{\Z})$ is zero.}

{\bf Proof.} Let us denote the ring $\Lambda_{\Z} \left< x_1,\ldots,x_g,x'_1,\ldots,x'_g \right>\otimes \Z[y] $ by $E^*_{g,n}$. Fix arbitrary $g\ge 1, n\ge 2$. Our aim is to prove that $\Tor(E^s_{g,n}/I^s_{Mac}) = 0, \forall s\ge 0$. \\
For any subset $\{i_1,\ldots,i_a,j_1,\ldots,j_b,k_1,\ldots,k_c\}\subset \{1,2,\ldots,g\}$ and any $q\ge 0$ let us consider the following two elements of $E^*_{g,n}$:
$$
\mathrm{monomial}  \  P = x_{i_1}\ldots x_{i_a}x'_{j_1}\ldots x'_{j_b}x_{k_1}x'_{k_1}\ldots x_{k_c}x'_{k_c}y^q,
$$
$$
\mathrm{and} \ S(P) = x_{i_1}\ldots x_{i_a}x'_{j_1}\ldots x'_{j_b}(y - x_{k_1}x'_{k_1})\ldots (y - x_{k_c}x'_{k_c})y^q.
$$
Clearly, $\deg P=\deg S(P) = a+b+2c+2q$. Set $w(P):=a+b+2c+q =$ the ``weight'' of $P$.

The ideal $I^*_{Mac}\subset E^*_{g,n}$ is generated by elements $S(P)$ for all monomials $P$ with $w(P)=n+1$. For such monomials $P$ we have $\deg S(P) = w(P) + q = n + 1 + q \ge n+1$. It follows that $I^s_{Mac}=0, 0\le s\le n.$ Therefore, our ring $E^*_{g,n}/I^*_{Mac}$ is torsion-free in dimensions $0\le s\le n$. 

{\bf Case $0\le s\le n$}. It is just done. 

{\bf Case $n+1\le s\le 2n-1$}. 

As before, we denote by $B_t, 0\le t\le 2n,$ the t-th Betti number of $\Sym^n M^2_g$. In Macdonald's paper there is an explicit formula for $B_t$:
$$
B_{2n-t} = B_t = \binom{2g}{t} + \binom{2g}{t-2} + \binom{2g}{t-4} + \ldots, \ \forall 0\le t\le n.
$$

First, let us prove that elements $S(P)$ lie in $I^*_{Mac}$ for all monomials $P$ with $w(P)\ge n+1$. We will use the induction on $w(P)$. The base $w(P)=n+1$ is correct.

The inductive step. Suppose that the statement is proved for $w(P)=n+t$ for some $t\ge 1$. Let us prove the statement for $w(P)=n+t+1$. Take an arbitrary monomial 
$$
P  = x_{i_1}\ldots x_{i_a}x'_{j_1}\ldots x'_{j_b}x_{k_1}x'_{k_1}\ldots x_{k_c}x'_{k_c}y^q,
$$
with $w(P)=a+b+2c+q =n+t+1$.

{\bf (1) $c\ge 1$}. In this case the corresponding element $S(P)$ has the following decomposition:
$$
S(P) = x_{i_1}\ldots x_{i_a}x'_{j_1}\ldots x'_{j_b}(y - x_{k_1}x'_{k_1})\ldots (y - x_{k_c}x'_{k_c})y^q = 
$$
$$
x_{i_1}\ldots x_{i_a}x'_{j_1}\ldots x'_{j_b}(y - x_{k_1}x'_{k_1})\ldots (y - x_{k_{c-1}}x'_{k_{c-1}})y^{q+1} \  - 
$$
$$
-  \  x_{k_c}x'_{k_c}x_{i_1}\ldots x_{i_a}x'_{j_1}\ldots x'_{j_b}(y - x_{k_1}x'_{k_1})\ldots (y - x_{k_{c-1}}x'_{k_{c-1}})y^q \  = \  S(P_1) + x'_{k_c}S(P_2),
$$
where 
$$
P_1 = x_{i_1}\ldots x_{i_a}x'_{j_1}\ldots x'_{j_b}x_{k_1}x'_{k_1}\ldots x_{k_{c-1}}x'_{k_{c-1}}y^{q+1} \  \mathrm{and} 
$$
$$
P_2 = x_{k_c}x_{i_1}\ldots x_{i_a}x'_{j_1}\ldots x'_{j_b}x_{k_1}x'_{k_1}\ldots x_{k_{c-1}}x'_{k_{c-1}}y^q.
$$
Evidently, $w(P_1)=w(P_2)=w(P)-1=n+t$. By induction hypothesis, we have $S(P_1),S(P_2)\in I^*_{Mac}$. It follows that $S(P)$ also lies in $I^*_{Mac}$. The inductive step in the case (1)  is proved. 

{\bf (2.1) $c=0,a=b=0$}. In this case $P=y^q=S(P)$ and $w(P)=q=n+t+1\ge n+1$. But, one has $y^{n+1}\in I^*_{Mac}$. Therefore, $S(P)$ also belongs to $I^*_{Mac}$. The inductive step in this case is proved. 

{\bf (2.2) $c=0,a\ge 1$}. In this case the corresponding element $S(P)$ has the following simple decomposition:
$$
S(P) = \pm x_{i_a} S(Q), 
$$
where
$$
Q =  x_{i_1}\ldots x_{i_{a-1}}x'_{j_1}\ldots x'_{j_b}y^q.
$$
Clearly, $w(Q)=w(P) - 1=n+t$. So, one has $S(Q)\in I^*_{Mac}$ and $S(P)\in I^*_{Mac}$. The inductive step in this case is proved. 

{\bf (2.3) $c=0,b\ge 1$}. The reasoning here is an evident change of that in the previous case.  

So, we have just proved the induction step in all cases. 

Therefore, we have proved, that $S(P)\in I^*_{Mac}$ for all monomials $P$ with $w(P)\ge n+1$.

Let us call a monomial $P$ ``primitive'', if $1\le w(P)\le n$. How many are primitive monomials $P$ of degree $s=n+t, 1\le t\le n-1$? We have two conditions for $P$:
$$
1\le w(P)=a+b+2c+q\le n,
$$
$$
\deg P=w(P)+q=n+t.
$$
It is clear, that for such $P$ we have $t\le q\le n+t -1$. Moreover, for $q=t$ there are $\binom{2g}{n-t}$ possibilities, for $q=t+1$ there are $\binom{2g}{n-t-2}$ possibilities, e.t.c. Thus, the whole quantity of primitive monomials $P$ of degree $n+t, 1\le t\le n-1,$ is equal to $\binom{2g}{n-t} + \binom{2g}{n-t-2} + \binom{2g}{n-t-4} +\ldots$. But, this number is equal to Betti number $B_{n+t}=B_{n-t}$.

Consider any monomial $P$ of degree $n+t, 1\le t\le n-1,$ and satisfying the condition $w(P)\ge n+1$. We have $S(P)\in I^*_{Mac}$. There are two possibilities: $(A)$ if $c=0$, then $P=S(P)\in I^*_{Mac}$; $(B)$ if $c\ge 1$, then by opening all brackets in $S(P)$ we obtain
$$
I^*_{Mac}\ni S(P) = \pm P \ + \  \sum_{ \mathrm{fin. \  num. \  of \ monomials \  with \ } a,b,c-1,q+1} \pm Q^{(1)}_{*} \ +
$$
$$
+ \ \sum_{ \mathrm{fin. \  num. \  of \ monomials \  with \ } a,b,c-2,q+2} \pm Q^{(2)}_{*} \ + \ldots.
$$
One has $w(Q^{(1)}_*)=w(P)-1, w(Q^{(2)}_*)=w(P)-2,$ e.t.c. By analysing the above identity, we obtain that in this case $(c\ge 1)$ the following representation holds:
$$
P = \sum_{\mathrm{fin.}} m_*Q_*  + R,
$$  
where $m_*\in \Z, \  Q_*$ are monomials of degree $n+t$, $w(Q_*)\le w(P) -1$, and $R\in I^*_{Mac}$.

By finite descending induction, for any monomial $P$ of degree $n+t, 1\le t\le n-1,$ and $w(P)\ge n+1$  we get the following decomposition:
$$
P = \sum_{\mathrm{fin.}} m_*Q_*  + R,
$$
where $m_*\in \Z, \  Q_*$ are monomials of degree $n+t$, $1\le w(Q_*)\le n$, and $R\in I^*_{Mac}$.

Therefore, for $1\le t\le n-1$ we obtain the identity for abelian groups:
$$
E^{n+t}_{g,n}/ I^{n+t}_{Mac} = \Z\left< Q | Q \mbox{ are primitive monomials of degree } n+t\right>/ (\Z\left< Q | \ldots \right>\cap I^{n+t}_{Mac}). 
$$
But, the free abelian group $\Z\left< Q | \ldots \right>$ has the rank equal to Betti number 
$$
B_{n+t} = \dim_{\Q}[(E^{n+t}_{g,n}/ I^{n+t}_{Mac})\otimes \Q].
$$
From the above two identities we get that $\Z\left< Q | \ldots \right>\cap I^{n+t}_{Mac} = 0$. Therefore, we have just proved that in the case $s=n+t,1\le t\le n-1,$ our abelian group $E^{s}_{g,n}/ I^{s}_{Mac}$ is torsion-free.

{\bf Case $s\ge 2n+1$}. Here we need to show that $E^s_{g,n}=I^s_{Mac}$.

Consider any monomial $P$ of degree $s\ge 2n+1$. We have $s=w(P) + q\ge 2n+1$. So, $w(P)\ge n+1 + (n-q)$ and $w(P)\ge q$. Thus, we always have $w(P)\ge n+1$ and $S(P)\in I^s_{Mac}$. 

Now, we also have two possibilities: $(A) \  c=0$, here $P=S(P)\in I^s_{Mac}$; and $(B) \ c\ge 1$, where we have the same decomposition of $S(P)$ as in the previous case $(n+1\le s\le 2n-1)$. 

Thus, for any monomial $P$ of degree $s\ge 2n+1$ we obtain $P\in I^s_{Mac}$, or we get the following decomposition:
$$
P = \sum_{\mathrm{fin.}} m_*Q_*  + R,
$$
where $m_*\in \Z, \  Q_*$ are monomials of degree $s$, $n+1\le w(Q_*)\le w(P) - 1$, and $R\in I^s_{Mac}$. 

By evident finite descending induction we obtain that any monomial $P$ of degree $s\ge 2n+1$ belongs to $I^s_{Mac}$. Therefore, we have just proved the desired identity $E^s_{g,n}=I^s_{Mac}$ in the considered case $(s\ge 2n+1)$.

{\bf Case $s=2n$}.  Let us show, that in this case there is only one primitive monomial $P=y^n$. Consider any primitive $P$ of the type $(a,b,c,q)$. We have two conditions:
$$
\deg P = a+b+2c +2q = 2n,
$$
$$
w(P) = a+b+2c +q \le n.
$$
Clearly, we get $0\le q\le n.$ If $q=n$, then $a=b=c=0$ and $P=y^n$. If $0\le q\le n-1$, then $w(P)= 2n - q\ge n+1$. Thus, we have just proved that the only primitive monomial in our last case is $P=y^n$.

Consider any monomial $P$ of the type $(a,b,c,q)$ and $\deg P=2n$. If $P\ne y^n$, then $w(P)\ge n+1$. Again we have two cases: $(A) \  c=0, P=S(P)\in I^*_{Mac}$; $(B) \ c\ge 1,$ here we obtain the following decomposition:
$$
P = \pm S(P) + \sum_{\mathrm{fin.}} m_*Q_*,
$$
where $m_*\in \Z, \  Q_*$ are monomials of degree $2n$, $n\le w(Q_*)\le w(P) - 1$. 

By evident finite descending induction we obtain that any monomial $P$ of degree $s=2n$ has the following decomposition:
$$
P = my^n + R,
$$
where $m\in \Z$ and $R\in I^*_{Mac}$.

Thus, we have the identity: 
$$
E^{2n}_{g,n}/ I^{2n}_{Mac} = \Z\left< y^n \right>/  (\Z\left< y^n \right>\cap I^{2n}_{Mac}).
$$ 
But, we know that $B_{2n} = \dim_{\Q} [(E^{2n}_{g,n}/ I^{2n}_{Mac}) \otimes \Q] = 1$. It follows that $\Z\left< y^n \right>\cap I^{2n}_{Mac} = 0$. Therefore, we have just proved that in the last case $(s=2n)$ our abelian group $E^{2n}_{g,n}/ I^{2n}_{Mac}$ is torsion-free. The proposition is completely proved.  \ \ \ \ \ \ \ \ \ \ \ \ \ \ \ \ \ \ \ \ \ \ \ \ \ \ \ \ \ \ \ \ \  \ \ \ \ \ \  \ \ \ \ \ \ \ \ \ \ \ \ \ \ \ \ \ \ \ \ \ \ \ \ \ \ \ \ \ \ \ \ \ \ \ \ \ \ \ \ \ \ \ \ \ \ \ \ \ \ \ \ \ \ \ \ \ \ \ \ \ \ \ \ \ \ \ \ \ \ \ \   \ \ \  $\Box$

\bigskip

\centerline{ \large ACKNOWLEDGEMENTS}   

\bigskip

The author is deeply grateful to his Advisor V.M.Buchstaber, A.A.Gaifullin, T.E.Panov and N.Ray for fruitful discussions.  The author is grateful to M.V.Prasolov for pointing out  D.Sullivan's theorem from \cite{Sul} and to S.Kallel and A.Polishchuk for pointing out R.Seroul's paper \cite{Ser}. The author is also deeply grateful to both referees for lots of constructive critics.

\begin{flushleft}
{\it 
Department of Geometry and Topology, \\
Faculty of Mechanics and Mathematics, \\
Moscow State University, Moscow, Russia \\
E-mail: dmitry-gugnin@yandex.ru
}
\end{flushleft}

\end{document}